\documentclass[11pt]{article}

\usepackage[T1]{fontenc}
\usepackage{lmodern}
\usepackage{microtype}
\usepackage[a4paper,margin=2.3cm]{geometry}
\usepackage{amsmath,amssymb,amsthm,mathtools}
\usepackage{enumitem}
\usepackage{xcolor}
\usepackage[colorlinks=true,linkcolor=blue!55!black,
 citecolor=blue!55!black,urlcolor=blue!55!black]{hyperref}
\usepackage[nameinlink,capitalize,noabbrev]{cleveref}
\hypersetup{
 pdftitle={Entropy Expansion for General Polynomial Images of Frostman Random Variables},
 pdfauthor={Alex Iosevich, Thang Pham, Nguyen Dac Quan, Steven Senger, and Boqing Xue},
 pdfsubject={Dyadic entropy expansion under Frostman hypotheses},
 pdfkeywords={Frostman measures, entropy, polynomial images, Elekes--Ronyai theorem, Fourier integral operators}
}

\newtheorem{theorem}{Theorem}[section]
\newtheorem{proposition}[theorem]{Proposition}
\newtheorem{lemma}[theorem]{Lemma}
\newtheorem{corollary}[theorem]{Corollary}
\newtheorem{definition}[theorem]{Definition}
\newtheorem{remark}[theorem]{Remark}
\newtheorem{example}[theorem]{Example}

\title{Entropy Expansion for General Polynomial Images\\
of Frostman Random Variables}
\author{Alex Iosevich \quad Thang Pham \quad Nguyen Dac Quan\\[0.5ex]
Steven Senger \quad Boqing Xue}
\date{}

\begin{document}
\maketitle

\begin{abstract}
We prove dyadic entropy expansion for the observables \(X+Y\) and
\(f(X,Y)\) under Frostman nonconcentration hypotheses on a prescribed,
possibly dependent law.  For an integer \(n\ge1\), write
\(H_n(Z)=H(\lfloor 2^nZ\rfloor)\) for the base-two Shannon entropy
at resolution \(2^{-n}\); thus \(n\) indexes the fineness of the
dyadic discretization.  For every \(0<s<1\) and every real bivariate
polynomial outside the family \(H(ax+by+c)\), conditional Frostman
regularity with common exponent \(s\) forces one of the two output entropies to exceed \(sn\) by
a fixed positive multiple of \(n\).  Under rectangular Frostman
regularity, the same conclusion holds after including the average
marginal entropy in the maximum.  The proof combines the discretized
Elekes--R\'onyai and web theorems of Raz and Zahl with simultaneous
information profiles and an estimate for mass near critical layers.
A second argument, based on Pham's folding estimate and a mutual
information penalty, gives explicit gains.  If the conditional laws
of \(X\) given \(Y\) and of \(Y\) given \(X\) have Frostman
exponents \(s_1,s_2\in(0,1)\), respectively, with
\(s_1+s_2>4/3\), we obtain
\(\max\{H_n(X+Y),H_n(f(X,Y))\}\ge
(\frac{s_1+s_2}{2}+\varepsilon)n-O(1)\)
for an explicit \(\varepsilon>0\).  Here the baseline averages the
two Frostman exponents.  We also classify the exceptional
coordinate representations and construct obstructions for algebraic
affine directions with Frostman constants uniform in the scale.
\end{abstract}


\tableofcontents

\section{Introduction, definitions, and main results}
\label{sec:main}

\subsection{Background and motivation}

A recurring principle in additive combinatorics and geometric measure
theory is that a nonconcentrated set or measure cannot remain small
under two sufficiently transverse algebraic observables.  At a fixed
scale, Frostman nonconcentration is the measure-theoretic counterpart
of the discretized nonconcentration conditions used by Katz and Tao
\cite{KatzTao}.  Bourgain's discretized sum--product and projection
theorem \cite{BourgainProjection} established a central expansion
principle for such sets.  Nonlinear projection variants were later
developed by Shmerkin \cite{ShmerkinNonlinear} and, in the
polynomial and web setting used here, by Raz and Zahl \cite{RZ}.
The Elekes--R\'onyai dichotomy \cite{ER} identifies the rigid
additive and multiplicative polynomial forms that underlie the
exceptional cases.

Dyadic entropy provides a scale-sensitive language for the same
phenomenon.  Local entropy averages have proved effective in the
study of projections and differentiable images of fractal measures
\cite{HochmanShmerkin}, while inverse theorems for entropy growth
under convolution form another influential part of this circle of
ideas \cite{HochmanEntropy}.  More recently, M\'ath\'e and O'Regan
\cite{MatheORegan} used Shannon-type information inequalities to
obtain quantitative discretized sum--product estimates.  These works
provide methodological context rather than black-box inputs for the
proofs below.

The present paper develops the entropy framework of our earlier work
\cite{IPQSX}.  Theorems~1.2 and~1.3 of that work treat, respectively,
jointly Frostman laws for diagonal polynomials of a common degree,
and conditionally Frostman laws for polynomials obtained from such
diagonal forms by an invertible rational linear change of variables.
Here the polynomial may be arbitrary outside the affine-univariate
family \(H(ax+by+c)\), with no rationality restriction on its
coefficients.  The full interval \(0<s<1\) is retained, and the
explicit estimates below also allow unequal conditional exponents.

We use two complementary nonlinear-projection inputs.  The theorems
of Raz and Zahl \cite{RZ} give qualitative expansion throughout
\(0<s<1\), while Minh-Quy Pham's folding Fourier integral operator
estimate \cite{PhamFIO} gives explicit gains in a smaller exponent
range.  The work here is to convert their set-theoretic and
product-measure estimates into uniform entropy bounds for the
prescribed law.  The simultaneous information-profile reduction, the
critical-layer estimate, the passage to dependent couplings, and the
classification of the original and two sum-coordinate
representations provide this conversion.  The families
\(\mathfrak E_S\) and \(\mathfrak E_3\) introduced below are the exact
exceptional families for these coordinate tests.  We do not assert
that they are the exact families for which entropy expansion itself
can fail.

\subsection{Definitions and main results}

We write \(\log=\log_2\) and use base-two Shannon entropy.
The symbol \(\ln\) denotes the natural logarithm used in the
geometric coordinate and curvature formulas.  For
\(z\in\mathbb R^k\), let
\[
 D_nz=2^{-n}\lfloor 2^nz\rfloor
\]
coordinatewise, and write
\[
 H_n(Z)=H(D_nZ),\qquad
 \operatorname{I}_n(X;Y)=H_n(X)+H_n(Y)-H_n(X,Y).
\]
Here \(n\ge1\) is the dyadic resolution level: \(D_n\) rounds
each coordinate down to the grid of mesh \(2^{-n}\), and \(H_n\)
is measured in bits.
All random variables below have bounded range.  Constants implicit in
\(O(1)\) may depend on this fixed range, but never on the dyadic scale
\(n\).

\begin{definition}[Frostman conditions]
Let \(0<s,s_1,s_2<1\).  A real random variable \(X\) is
\((s,C)\)-Frostman if
\[
 \mathbb P(X\in I)\le C|I|^s
\]
for every interval \(I\) of length at most one.  A pair \((X,Y)\) is
jointly \((s_1,s_2,C)\)-Frostman if
\[
 \mathbb P(X\in I,Y\in J)
 \le C|I|^{s_1}|J|^{s_2}
\]
for all such intervals.  It is conditionally
\((s_1,C_1;s_2,C_2)\)-Frostman if, for all intervals \(I,J\) of
length at most one and all Borel sets \(E\),
\begin{align*}
 \mathbb P(X\in I,Y\in E)&\le C_1|I|^{s_1}\mathbb P(Y\in E),\\
 \mathbb P(Y\in J,X\in E)&\le C_2|J|^{s_2}\mathbb P(X\in E).
\end{align*}
\end{definition}

The conditional hypothesis implies the joint one, with constant
\(C_1C_2\).  It also implies the corresponding marginal Frostman
bounds.  We shall repeatedly use the following equivalent
disintegrated consequence.  For almost every \(y\), the conditional
law of \(X\) given \(Y=y\) is \((s_1,C_1)\)-Frostman, and for almost
every \(x\), the conditional law of \(Y\) given \(X=x\) is
\((s_2,C_2)\)-Frostman.  Indeed, for each interval with rational
endpoints the defining event inequality gives the corresponding
Radon--Nikodym bound; intersecting these countably many full-measure
sets and approximating an arbitrary interval proves the assertion.

\medskip
\noindent{Geometric comparison.
Let \(T\) be uniformly distributed on \([0,1]\) and put
\((X,Y)=(T,T)\).  Then, for intervals \(I,J\) of length at most one,
\[
 \mathbb P(X\in I,Y\in J)
 =|I\cap J\cap[0,1]|
 \le \min\{|I|,|J|\}
 \le |I|^{1/2}|J|^{1/2}.
\]
Thus \((X,Y)\) is jointly \((1/2,1/2,1)\)-Frostman.  It is not
conditionally Frostman with any positive exponent: for almost every
\(x\), the conditional law of \(Y\) given \(X=x\) is the Dirac mass
\(\delta_x\).  Hence joint Frostman regularity permits concentration
on a non-axis curve, whereas conditional Frostman regularity rules
out this fiberwise collapse.}

Following Elekes--R\'onyai \cite{ER}, we call a polynomial
\(f\in\mathbb R[x,y]\) a
polynomial special form if
\begin{equation}
 f(x,y)=H(A(x)+B(y))
 \quad\text{or}\quad
 f(x,y)=H(A(x)B(y))
 \label{eq:polynomial-special}
\end{equation}
for univariate real polynomials \(H,A,B\).  A real-analytic function
\(F\) on a connected open set \(\Omega\) is
analytic-special in the Raz--Zahl--Pham sense if, on every
connected component of
\(
 \Omega\setminus(\{F_x=0\}\cup\{F_y=0\}),
\)
it can be written
\begin{equation}
 F(x,y)=g(h(x)+k(y))
 \label{eq:analytic-special}
\end{equation}
with one-variable real-analytic functions \(g,h,k\).  This definition
includes multiplicative forms after taking natural logarithms on a
regular sign component.  For polynomials the analytic and polynomial
notions agree.  If none of \(f_x,f_y,f_{xy}\) vanishes identically,
this is the polynomial curvature criterion in
\cite[Lemma~3.3]{RZ}, together with its analytic version.  If \(f_x\)
or \(f_y\) vanishes identically, the polynomial is univariate.  If
\(f_{xy}\equiv0\), it has the additive form \(A(x)+B(y)\).
These cases establish the equivalence without any nondegeneracy
qualification.  We retain both terms to match the formulations of
the external inputs and to discuss nonpolynomial functions in local
charts.  For a polynomial on \(\mathbb R^2\), analytic-specialness always
refers to its regular components in \(\mathbb R^2\).  The separated
derivative-ratio identity below and analytic continuation show that
one nonempty regular product rectangle suffices to verify it; when
the regular set is empty, the univariate case applies.

We shall use the two exceptional families
\begin{align}
 \mathfrak E_S
 &=\{H(ax+P(x+y)):H,P\in\mathbb R[t],\ a\in\mathbb R\},
 \label{eq:ES-main}\\
 \mathfrak E_3
 &=\{H(ax+by+c):H\in\mathbb R[t],\ a,b,c\in\mathbb R\}.
 \label{eq:E3-main}
\end{align}
Thus \(\mathfrak E_3\subset\mathfrak E_S\).  The subscript in
\(\mathfrak E_S\) refers to the sum coordinate \(S=x+y\), whereas
\(\mathfrak E_3\) records the simultaneous use of the original and
the two sum-coordinate systems.

\medskip
\noindent{Geometric meaning.
A nonconstant polynomial in \(\mathfrak E_3\) factors through one
affine coordinate.  Consequently, each of its level sets is empty
or a finite union of parallel affine lines; for example,
\(f(x,y)=(x+2y)^2\) lies in \(\mathfrak E_3\).  The larger family
\(\mathfrak E_S\) is adapted to the sum coordinate \(S=x+y\): if
\(f(x,y)=H(ax+P(S))\), then
\[
 G_X(x,S)=H(ax+P(S)),\qquad
 G_Y(y,S)=H(-ay+P(S)+aS),
\]
so both sum transitions remain analytic-special.  The classification
results below show that, among nonconstant polynomials,
\(\mathfrak E_S\) is exactly the family for which both transitions are
analytic-special, while \(\mathfrak E_3\) is exactly the family for
which the original map and both transitions are all analytic-special.}

Our first theorem is qualitative but holds throughout the full
Frostman range.

\begin{theorem}[General-polynomial entropy expansion]
\label{thm:RZ-general}
Let \(f\in\mathbb R[x,y]\setminus\mathfrak E_3\), and fix
\(0<s<1\) and \(C,c>0\).
\begin{enumerate}[label=\textup{(\roman*)}]
\item There are \(\varepsilon=\varepsilon(f,s,c)>0\) and
\(N=N(f,s,C,c)\) such that every conditionally
\((s,C;s,C)\)-Frostman pair \((X,Y)\), supported in
\([-c,c]^2\), satisfies
\begin{equation}
 \max\{H_n(X+Y),H_n(f(X,Y))\}\ge(s+\varepsilon)n
 \label{eq:RZ-general-cond}
\end{equation}
for every \(n>N\).

\item There are \(\varepsilon=\varepsilon(f,s,c)>0\) and
\(N=N(f,s,C,c)\) such that every jointly
\((s,s,C)\)-Frostman pair \((X,Y)\), supported in
\([-c,c]^2\), satisfies
\begin{equation}
 \max\left\{\frac{H_n(X)+H_n(Y)}{2},
 H_n(X+Y),H_n(f(X,Y))\right\}\ge(s+\varepsilon)n
 \label{eq:RZ-general-joint}
\end{equation}
for every \(n>N\).
\end{enumerate}
\end{theorem}

The marginal-entropy term in \eqref{eq:RZ-general-joint} is essential
in general.  Let \(T\) be uniform on \([0,1]\), put
\((X,Y)=(T,-T)\), and take \(f(x,y)=x^2-y^2\).  For intervals of
length at most one,
\[
 \mathbb P(X\in I,Y\in J)
 =|I\cap(-J)\cap[0,1]|
 \le\min\{|I|,|J|\}
 \le |I|^{1/2}|J|^{1/2}.
\]
Thus the law is jointly \((1/2,1/2,1)\)-Frostman.  The polynomial
lies outside \(\mathfrak E_3\): its gradient \((2x,-2y)\) has no
nonzero constant annihilating direction, as characterized in
\cref{cor:E3-test}.  Nevertheless,
\[
 X+Y=f(X,Y)=0,
 \qquad H_n(X)=H_n(Y)=n.
\]
The two output entropies therefore vanish, while the marginal term
supplies the expansion.  This example also distinguishes the joint
conclusion from the conditional one.

For \(s>1/2\), the jointly Frostman assertion can be obtained from
the localized curved three-web and flat four-web theorems alone:
the \(r\)-neighborhood of a fixed real algebraic curve has mass
\(O(r^{2s-1})\), so a sufficiently thin fixed tubular neighborhood
may be discarded with arbitrarily small mass.  For arbitrary
\(0<s<1\), such deletion is not available for a merely jointly
Frostman law; in the branch where the polynomial is not a polynomial
special form, the proof instead uses Raz--Zahl's global discretized
Elekes--R\'onyai theorem.  Under conditional Frostman regularity,
algebraic tube deletion remains available for every positive
exponent.  Raz--Zahl's expansion exponent is polynomial-independent,
but the admissible density tolerance in their global theorem is stated
with dependence on the polynomial.  The safe conclusion is therefore
the displayed \(\varepsilon(f,s,c)\), rather than a degree-uniform constant;
the dependence on \(c\) records the affine normalization of the support.

The second theorem records the explicit gains obtained from Pham's
folding estimate.  The symbol \(O(1)\) is uniform in \(n\) and in the
admissible law; it may depend on the polynomial, the displayed
exponents, the Frostman constants, the range, and the chosen strict
entropy exponent.

\begin{theorem}[General-polynomial expansion with explicit gains]
\label{thm:Pham-explicit}
Let \(f\in\mathbb R[x,y]\), fix \(c>0\), and suppose all laws below are
supported in \([-c,c]^2\).
\begin{enumerate}[label=\textup{(\roman*)}]
\item If \((X,Y)\) is jointly \((s,s,C)\)-Frostman, \(f\) is not
analytic-special, and \(2/3<s<1\), then
\[
 \max\left\{\frac{H_n(X)+H_n(Y)}{2},H_n(f(X,Y))\right\}
 \ge(s+\varepsilon)n-O(1)
\]
for every integer \(n\ge1\) and every
\[
 0<\varepsilon<\frac{1}{5}\min\left\{s-\frac{2}{3},1-s\right\}.
\]

\item If \((X,Y)\) is jointly \((s,s,C)\)-Frostman,
\(f\notin\mathfrak E_3\), and \(5/6<s<1\), then
\begin{equation}
 \max\left\{\frac{H_n(X)+H_n(Y)}{2},H_n(X+Y),
 H_n(f(X,Y))\right\}
 \ge(s+\varepsilon)n-O(1)
 \label{eq:Pham-joint-main}
\end{equation}
for every integer \(n\ge1\) and every
\[
 0<\varepsilon<\frac{1}{7}
 \min\left\{2s-\frac{5}{3},1-s\right\}.
\]
In particular, the safe explicit value
\begin{equation}
 \varepsilon=\frac{1}{14}
 \min\left\{2s-\frac{5}{3},1-s\right\}.
 \label{eq:Pham-joint-eps}
\end{equation}
is valid.  If both sum transitions in
\eqref{eq:sum-transitions} below are not analytic-special, the full
range improves by replacing \(1/7\) with \(1/5\); correspondingly,
the safe factor \(1/14\) may be replaced by \(1/10\).

\item Let \((X,Y)\) be conditionally
\((s_1,C_1;s_2,C_2)\)-Frostman, where \(0<s_1,s_2<1\).  Put
\[
 \bar s=\frac{s_1+s_2}{2},\qquad
 s_-=\min\{s_1,s_2\},\qquad
 q=\min\left\{1,s_1+s_2-\frac{2}{3}\right\}.
\]
If \(s_1+s_2>4/3\) and \(f\notin\mathfrak E_3\), then
\begin{equation}
 \max\{H_n(X+Y),H_n(f(X,Y))\}
 \ge(\bar s+\varepsilon)n-O(1)
 \label{eq:Pham-cond-main}
\end{equation}
for every integer \(n\ge1\) and every
\[
 0<\varepsilon<\frac{q-\bar s}{5}
 =\frac{1}{5}\min\left\{\bar s-\frac{2}{3},1-\bar s\right\}.
\]
In particular, \(\varepsilon=(q-\bar s)/10\) is valid.  If, more strongly,
\(f\notin\mathfrak E_S\), one also has the complementary estimate
\begin{equation}
 \max\{H_n(X+Y),H_n(f(X,Y))\}
 \ge(s_-+\varepsilon)n-O(1)
 \label{eq:Pham-cond-transition-main}
\end{equation}
for every \(0<\varepsilon<(q-s_-)/3\); in particular,
\(\varepsilon=(q-s_-)/6\) is valid in this estimate.  More precisely, if
the transition \(G_X\) in \eqref{eq:sum-transitions} is not
analytic-special, then
\begin{equation}
 \max\{H_n(X+Y),H_n(f(X,Y))\}
 \ge\max\left\{s_1,\frac{t+2s_2}{3}\right\}n-O(1)
 \label{eq:Pham-cond-oriented}
\end{equation}
for every \(0<t<q\); the symmetric statement holds for \(G_Y\),
with \(s_1,s_2\) interchanged.  These oriented bounds imply both
\eqref{eq:Pham-cond-main} and
\eqref{eq:Pham-cond-transition-main} in the transition case.
When \(s_1=s_2=s\), the assumptions reduce to \(2/3<s<1\), and the
uniform choice for \(f\notin\mathfrak E_3\) is
\begin{equation}
 \varepsilon=\frac{1}{10}\min\left\{s-\frac{2}{3},1-s\right\}.
 \label{eq:Pham-cond-eps}
\end{equation}
The factor \(1/10\) may be replaced by \(1/6\) when
\(f\notin\mathfrak E_S\).

\item If \(X,Y\) are independent, \(X\) is
\((s_1,C_1)\)-Frostman, \(Y\) is \((s_2,C_2)\)-Frostman,
\(0<s_1,s_2<1\), and \(f\) is not analytic-special, then
\[
 H_n(f(X,Y))\ge u n-O(1)
\]
for every integer \(n\ge1\) and every
\(0<u<\min\{1,s_1+s_2-2/3\}\).  If instead only
\(f\notin\mathfrak E_3\) is assumed and \(s_1+s_2>4/3\), then
\eqref{eq:Pham-cond-main} holds with \(\bar s\) and \(q\) as in
\textup{(iii)}, for every integer \(n\ge1\), and with the same
range of \(\varepsilon\).
\end{enumerate}
\end{theorem}

The two theorems answer different quantitative questions.
\Cref{thm:RZ-general} reaches every \(0<s<1\) but gives no useful
closed formula for the gain.  \Cref{thm:Pham-explicit} supplies explicit
gains, with the combined cross-partition/fold threshold loss
\(2/3=(2-1)/2+1/6\); the fold derivative loss itself is \(1/6\).
The remaining sections state the external inputs precisely and prove
the reductions leading to these conclusions.

\subsection{Examples covered by the general theorem}

The polynomial scope can be seen without evaluating a curvature
expression.  By \cref{cor:E3-test}, it is enough to check that
\(uf_x+vf_y\equiv0\) forces \(u=v=0\).  Thus the qualitative
theorem applies to \(x^2+y\), to \(x^d+y^d\) for every \(d\ge2\),
and to every real quadratic polynomial whose quadratic part is
nondegenerate.  No rationality condition is imposed on the quadratic
coefficients.  The exact quadratic criterion is given in
\cref{prop:diagonal-quadratic-examples}.

A mixed cubic illustrating the extension beyond the family in
\cite{IPQSX} is
\[
 f(x,y)=xy(x+y).
\]
Here \(f_x=2xy+y^2\) and \(f_y=x^2+2xy\), so comparison of the
\(x^2\) and \(y^2\) coefficients proves \(f\notin\mathfrak E_3\).
It cannot be written as a sum of two univariate polynomials in
independent real linear coordinates.  Indeed, its highest homogeneous
part has three distinct real linear factors.  In a representation
of that kind, both univariate degrees would have to be at most three.
If both cubic coefficients were nonzero, the highest part would be
a sum of two cubes of independent real linear forms, which has only
one real projective zero.  If one cubic coefficient vanished, the
highest part would be a single cube.  Both alternatives contradict
the factorization \(xy(x+y)\).

The polynomial \((x+y)^2+x-y\) exhibits a different feature: both
sum-coordinate transitions are special, while the original map is
not.  It shows why a treatment based only on the two sum transitions
is insufficient.  The calculations in
\cref{ex:original-coordinates-needed} explain this distinction.

\subsection{Overview of the proofs}

We first describe the Raz--Zahl route.  Suppose that one of the
conclusions of \cref{thm:RZ-general} fails at a large scale
\(\delta=2^{-n}\).  Under the conditional hypothesis,
\cref{thm:coordinate-entropy} recovers the source entropy from the
pair of outputs \((X+Y,f(X,Y))\); under the joint hypothesis, the
source and marginal upper bounds follow directly from the failed
maximum and subadditivity.  Frostman lower bounds then confine the
source and marginal information, on most of the probability space,
to short windows around their natural exponents.  The simultaneous
profile lemma, \cref{lem:information-profiles}, refines the source
cells by all relevant output cells, discards the information tails by
shifted Markov inequalities, and groups the remaining atoms according
to these windows.  Each retained profile produces a set
\(E_\pi\) with the required size and nonconcentration estimates, while
the covering numbers of all relevant images of \(E_\pi\) are controlled
by the same profile.  This simultaneity is the bridge from averaged
entropy bounds to the geometric hypotheses in the Raz--Zahl inputs.

If \(f\) is not a polynomial special form, the global discretized
Elekes--R\'onyai theorem \cref{prop:RZ-global} forces
\(f(E_\pi)\) to have a larger covering exponent.  Averaging the
resulting output-information lower bounds contradicts the assumed
entropy failure.  This argument works for jointly Frostman laws
throughout \(0<s<1\), since it does not remove a neighborhood of a
curved algebraic set.  The localized curved-web theorem
\cref{prop:RZ-curved} gives the alternative proof in
\cref{sec:curved-appendix} under conditional
Frostman regularity for every positive exponent and under joint
Frostman regularity when \(s>1/2\), precisely the range in which a
fixed algebraic tube has small mass.

For a polynomial special form outside \(\mathfrak E_3\), regular
product charts supplied by \cref{lem:RZ-special-charts} reduce the
two observables to
\[
 f(x,y)=h(p(x)+q(y)),
 \qquad
 x+y=F(p(x))+G(q(y)).
\]
The flat four-web theorem \cref{prop:RZ-flat} is applied
simultaneously to \(p,q,p+q\), and \(F(p)+G(q)\).  The derivative of
the outer function \(h\) may vanish, and a merely jointly Frostman law
may place substantial mass near the corresponding critical layers.
Rather than deleting those layers, \cref{lem:critical-layer} shows
that large critical-layer mass itself forces excess marginal entropy.
Away from the layers, the web gain passes from \(p+q\) to \(f\) with
a controlled loss.  A linear inequality combines the pairwise
coordinate bounds with the web gain, and averaging over the retained
profiles gives the required contradiction.

The explicit theorem starts from a different input.  Away from the
algebraic degeneracy set of a non-special polynomial,
\cref{prop_pham_local_energy} converts Pham's two-sided fold estimate
into
\[
 I_u\bigl(f_*(\mu_1\times\mu_2)\bigr)
 \lesssim
 \bigl(1+I_{a_1}(\mu_1)\bigr)
 \bigl(1+I_{a_2}(\mu_2)\bigr),
 \qquad a_1+a_2>u+\frac{2}{3}.
\]
This estimate concerns the product of the marginals, whereas the law
in the theorem may be dependent.  The relative-entropy inequality in
\cref{lem_pham_relative_entropy}, followed by data processing and a
bounded-displacement comparison, yields the correlation penalty
\cref{prop_pham_correlation}:
\[
 H_n(f(U,V))+2\operatorname{I}_{n+r}(U;V)\ge tn-O(1).
\]
Thus failure of output expansion forces large mutual information.
On the other hand, the Frostman lower bound for the joint source entropy,
together with the candidate upper bounds in the relevant maximum,
bounds the same mutual information from above.  Combining these two
bounds produces a linear inequality with an explicit entropy gain.

Finally, \cref{thm:two-transition,cor:three-coordinate} identify the
coordinate systems in which the FIO input is non-special.  Unless
\(f\in\mathfrak E_3\), at least one of the original representation
and the two sum transitions is non-special.  For a sum transition,
the invertible shear \((X,Y)\mapsto(X,X+Y)\), or its symmetric
counterpart, preserves the joint source entropy and expresses the
correlation penalty in terms of the three entropies in
\cref{thm:Pham-explicit}.  If both transitions are non-special, the two
correlation inequalities can be combined and the coefficient
improves.  Combining these linear inequalities gives the factors
\(1/5\), \(1/7\), and their conditional analogues.  In the
asymmetric conditional case, the additional bound
\(H_n(X+Y)\ge\max\{s_1,s_2\}n-O(1)\) converts the transition
estimate into a gain above the average exponent \(\bar s\).  This also explains
the complementary nature of the two methods: the Raz--Zahl argument
reaches the full Frostman range with a qualitative gain, whereas the
FIO argument retains the explicit \(2/3\)-loss and therefore gives
explicit values only in the stated higher-exponent ranges.

\section{Raz--Zahl inputs and sum-coordinate transitions}
\label{sec:RZ-inputs}

For a bounded set \(E\subset\mathbb R^k\), let \(\mathcal N_\delta(E)\) denote
the least number of axis-parallel \(\delta\)-cubes needed to cover
\(E\).  We freely replace a union of \(\delta\)-cells by its set of
cell centers, or conversely thicken a \(\delta\)-separated set by
\(\delta\)-cells.  These operations change all covering estimates by
fixed constants only.  We now record the three consequences of
Raz--Zahl that will be used later.

\subsection{The global expanding-polynomial theorem}

\begin{proposition}[Raz--Zahl, global form]
\label{prop:RZ-global}
Let \(0<s<1\), and let \(g\) be a real polynomial defined on a
neighborhood of \([0,1]^2\) which is not a polynomial special form
in the sense of \eqref{eq:polynomial-special}.  There are
\[
 \kappa_{\rm RZ}=\kappa_{\rm RZ}(s)>0,
 \qquad
 \eta_{\rm RZ}=\eta_{\rm RZ}(g,s)>0
\]
such that the following holds for all sufficiently small \(\delta\).
Let \(A,B\subset[0,1]\) be unions of \(\delta\)-intervals satisfying
\begin{equation}
 |A\cap J|,|B\cap J|
 \le |J|^s\delta^{1-s-\eta_{\rm RZ}}
 \label{eq:RZ-global-NC}
\end{equation}
for every interval \(J\).  If
\(E\subset A\times B\) is a union of \(\delta\)-squares and
\begin{equation}
 |E|\ge\delta^{2-2s+\eta_{\rm RZ}},
 \label{eq:RZ-global-size}
\end{equation}
then
\begin{equation}
 |g(E)|\ge\delta^{1-s-\kappa_{\rm RZ}}.
 \label{eq:RZ-global-conclusion}
\end{equation}
Here the two occurrences of \(|\cdot|\) in
\eqref{eq:RZ-global-size} and \eqref{eq:RZ-global-conclusion} denote
planar and one-dimensional Lebesgue measure, respectively.
\end{proposition}

\Cref{prop:RZ-global} is the specialization of
\cite[Theorem~1.18]{RZ} with its two source parameters equal to
\(s\).  Its expansion gain may be chosen independently of the fixed
polynomial, while the admissible nonconcentration tolerance and the
threshold scale may depend on the polynomial.  We therefore retain
polynomial dependence in the gains asserted in \Cref{thm:RZ-general}.

\subsection{The curved three-web theorem}

\begin{proposition}[Raz--Zahl, curved three-web form]
\label{prop:RZ-curved}
Fix \(0<s<1\), a compact set \(K\Subset U\subset\mathbb R^2\), and
smooth functions \(F_1,F_2,F_3:U\to\mathbb R\).  Suppose that their
gradients are pairwise transverse on \(U\) and that the Blaschke
curvature of the planar three-web
\((F_1,F_2,F_3)\) is nowhere zero on \(U\).  There are
\(\eta_{\rm curv},\kappa_{\rm curv}>0\) and \(\delta_0>0\),
where \(\eta_{\rm curv}\) and \(\kappa_{\rm curv}\) may be chosen
depending only on \(s\), such that the following holds.  The
threshold \(\delta_0\) may depend also on the fixed buffered chart.
If \(0<\delta<\delta_0\) and a union
\(E\subset K\) of \(\delta\)-squares satisfies
\begin{align}
 |E|&\ge\delta^{2-2s+\eta_{\rm curv}},
 \label{eq:RZ-curved-size}\\
 |E\cap B(z,r)|
 &\lesssim r^{2s}\delta^{2-2s-\eta_{\rm curv}}
 \quad(\delta\le r\le1),
 \label{eq:RZ-curved-NC}
\end{align}
then
\begin{equation}
 \max_{1\le i\le3}\mathcal N_\delta(F_i(E))
 \gtrsim\delta^{-s-\kappa_{\rm curv}}.
 \label{eq:RZ-curved-conclusion}
\end{equation}
Fixed multiplicative constants in the hypotheses are harmless after
slightly decreasing the two positive exponents.
\end{proposition}

This is the covering-number consequence of
\cite[Theorem~1.11]{RZ}, specialized to source exponent \(2s\).
Indeed, each image of a \(\delta\)-square under a fixed smooth map has
diameter \(O(\delta)\), and hence
\(|F(E)|\lesssim\delta\mathcal N_\delta(F(E))\).
The transversality and curvature bounds are uniform on the buffered
chart; the threshold scale may depend on those bounds and on the
three fixed functions.

\subsection{The flat four-web theorem}

\begin{proposition}[Raz--Zahl, flat four-web form]
\label{prop:RZ-flat}
Fix \(0<s<1\), compact intervals \(I,J\), and smooth functions
\(u:I\to\mathbb R\), \(v:J\to\mathbb R\).  Suppose that
\begin{equation}
 \inf_I|u'|>0,\qquad \inf_I|u''|>0,\qquad
 \inf_J|v'|>0.
 \label{eq:flat-regularity}
\end{equation}
There are \(\eta_{\rm web}>0\), \(\kappa_{\rm web}>0\), and
\(\delta_0>0\), where \(\eta_{\rm web}\) and
\(\kappa_{\rm web}\) may be chosen depending only on \(s\), such
that the following holds.  The threshold \(\delta_0\) may depend
also on the fixed buffered chart.  If \(0<\delta<\delta_0\) and a union
\(E\subset I\times J\) of \(\delta\)-squares satisfies
\begin{align}
 |E|&\ge\delta^{2-2s+\eta_{\rm web}},
 \label{eq:RZ-flat-size}\\
 |E\cap B(z,r)|
 &\lesssim r^{2s}\delta^{2-2s-\eta_{\rm web}}
 \quad(\delta\le r\le1),
 \label{eq:RZ-flat-NC}
\end{align}
then
\begin{equation}
 \max\left\{
 \mathcal N_\delta(x(E)),\mathcal N_\delta(y(E)),
 \mathcal N_\delta((x+y)(E)),
 \mathcal N_\delta((u(x)+v(y))(E))
 \right\}
 \gtrsim\delta^{-s-\kappa_{\rm web}}.
 \label{eq:RZ-flat-conclusion}
\end{equation}
\end{proposition}

This is the covering-number consequence of
\cite[Theorem~1.12]{RZ}.  In the technical flat theorem, it is
enough to bound \(|u'|,|u''|,|v'|\) away from zero; no lower bound on
\(|v''|\) is required.  The roles of the two coordinates may of
course be interchanged.  As above, the conversion from Lebesgue
measure to covering number uses
\(|F(E)|\lesssim\delta\mathcal N_\delta(F(E))\).

\begin{remark}[Buffered charts]
The word buffered means that the compact chart core is
contained in the interior of a slightly larger chart on which all
regularity and transversality bounds continue to hold.  For large
\(n\), every level-\(n\) cell meeting the core then lies in the
larger chart.  This removes boundary errors without changing any
entropy or covering exponent.
\end{remark}

\begin{lemma}[Buffered coordinate changes]
\label{lem:buffered-coordinate-change}
Let \(U_0\Subset U\subset\mathbb R^2\), and let
\(\Psi:U\to\Psi(U)\) be a fixed \(C^2\) diffeomorphism whose
derivative and inverse derivative are bounded on \(U\).  Put
\(\delta=2^{-n}\).  There is \(n_0\), depending only on the buffer,
such that, for every \(n\ge n_0\), every level-\(n\) square meeting
\(U_0\) is contained in \(U\).  If \(E\) is a union of such squares,
let \(\widetilde E\) be the union of the level-\(n\) squares contained
in \(\Psi(U)\) that meet \(\Psi(E)\).  Then, with constants depending
only on the buffered coordinate change,
\begin{align}
 |\widetilde E|&\asymp |E|,\label{eq:coordinate-change-area}\\
 |\widetilde E\cap B(z,r)|
 &\lesssim |E\cap B(z_0,Cr)|
 \qquad(\delta\le r\le1),
 \label{eq:coordinate-change-NC}
\end{align}
where, whenever the left-hand side is nonzero, one may choose
\(z_0\in U\) with \(|\Psi(z_0)-z|\lesssim r\).
For every fixed Lipschitz scalar map \(F\) on \(U\),
\begin{equation}
 \mathcal N_\delta\bigl((F\circ\Psi^{-1})(\widetilde E)\bigr)
 \asymp_F \mathcal N_\delta(F(E)).
 \label{eq:coordinate-change-output}
\end{equation}
Consequently, the hypotheses and conclusions of
\cref{prop:RZ-curved,prop:RZ-flat} are unchanged, up to fixed
constants and an arbitrarily small decrease of their positive
exponents, under a buffered coordinate change and a fixed affine
rescaling.
\end{lemma}

\begin{proof}
A level-\(n\) square has image contained in a ball of radius
\(O(\delta)\), and the same holds for \(\Psi^{-1}\).  Thus each source
square meets only \(O(1)\) squares after applying \(\Psi\), and each
target square receives points from only \(O(1)\) source squares.
This proves \eqref{eq:coordinate-change-area}.  The inverse image of
\(B(z,r+O(\delta))\) lies in a ball of radius \(Cr\), which gives
\eqref{eq:coordinate-change-NC}.  More explicitly, compactness of the
buffered core and the inverse-function theorem give this linear bound
uniformly for all sufficiently small \(r\); for larger \(r\), it
follows after increasing \(C\) by the fixed diameter of the core.
Finally, saturating either the
source or target by level-\(n\) squares changes a fixed Lipschitz
scalar image by an \(O(\delta)\)-neighborhood.  Bounded overlap in
both directions proves \eqref{eq:coordinate-change-output}.
\end{proof}

\subsection{The two sum-coordinate systems}

Put \(S=x+y\) and define the two invertible shears
\[
 T_X(x,y)=(x,S),\qquad T_Y(x,y)=(y,S).
\]
For \(f\in\mathbb R[x,y]\), its sum transitions are
\begin{equation}
 G_X(x,S)=f(x,S-x),
 \qquad
 G_Y(y,S)=f(S-y,y).
 \label{eq:sum-transitions}
\end{equation}
Fixed shears preserve dyadic entropy up to \(O(1)\), but they do not
preserve every Frostman hypothesis with the same exponents.  The
precise facts used later are the following.

\begin{lemma}[Conditional Frostman under a sum shear]
\label{lem:conditional-shear}
If \((X,Y)\) is conditionally
\((s_1,C_1;s_2,C_2)\)-Frostman and \(S=X+Y\), then
\((X,S)\) is jointly \((s_1,s_2,C')\)-Frostman and
\((Y,S)\) is jointly \((s_2,s_1,C')\)-Frostman, where
\(C'\) depends only on \(C_1,C_2\).  Moreover, for every \(m\),
\begin{align}
 H(D_mS\mid D_mX)&\ge s_2m-O(1),
 \label{eq:shear-cond-X}\\
 H(D_mS\mid D_mY)&\ge s_1m-O(1).
 \label{eq:shear-cond-Y}
\end{align}
\end{lemma}

\begin{proof}
For fixed \(x\), the condition \(x+Y\in J\) restricts \(Y\) to an
interval of length \(|J|\).  The conditional Frostman estimate for
\(Y\) given \(X\), followed by the marginal estimate for \(X\),
gives
\[
 \mathbb P(X\in I,S\in J)
 \le C_1C_2|I|^{s_1}|J|^{s_2}.
\]
The other shear is symmetric.  Given a level-\(m\) \(X\)-cell, the
dyadic ambiguity between \(S\) and \(Y=S-X\) is bounded.  The
conditional atom bound for \(Y\) therefore implies
\eqref{eq:shear-cond-X}; the other inequality is identical.
\end{proof}

\subsection{The algebra of the two transitions}

Recall the componentwise definition of analytic-specialness in
\eqref{eq:analytic-special}.  On a smaller regular product
neighborhood, analytic-specialness implies
\begin{equation}
 \frac{F_u}{F_v}=A(u)B(v).
 \label{eq:separated-derivative-ratio}
\end{equation}
Every polynomial special form \eqref{eq:polynomial-special} has this
property on its regular components.

\begin{theorem}[Two-transition classification]
\label{thm:two-transition}
Let \(f\in\mathbb R[x,y]\) be nonconstant, and let \(G_X,G_Y\) be as in
\eqref{eq:sum-transitions}.  Both \(G_X\) and \(G_Y\) are
analytic-special if and only if
\(f\in\mathfrak E_S\).
\end{theorem}

\begin{proof}
If \(f(x,y)=H(ax+P(x+y))\), then
\[
 G_X(x,S)=H(ax+P(S)),
 \qquad
 G_Y(y,S)=H(-ay+P(S)+aS),
\]
so both transitions are analytic-special.

Conversely, put \(p=f_x\), \(q=f_y\), and \(r=p/q\).  If one of
\(p,q,p-q\) vanishes identically, then \(f\) is a polynomial in
\(y,x\), or \(x+y\), respectively, and the conclusion is immediate.
Otherwise,
\[
 \frac{(G_X)_x}{(G_X)_S}=r-1,
 \qquad
 \frac{(G_Y)_y}{(G_Y)_S}=\frac{1-r}{r}.
\]
Equation \eqref{eq:separated-derivative-ratio}, first on regular
rectangles and then as a rational-function identity, supplies
nonzero one-variable rational functions \(A,B,C,D\) such that
\[
 r-1=A(x)B(S),
 \qquad
 \frac{1-r}{r}=C(y)D(S).
\]
To justify rationality of the separated factors, fix a regular base
point \((x_0,S_0)\).  If \(R(x,S)=r(x,S-x)-1\), normalize on a
smaller product rectangle by
\[
 A(x)=\frac{R(x,S_0)}{R(x_0,S_0)},
 \qquad B(S)=R(x_0,S).
\]
Both factors are rational, and the identity \(R=A B\), initially
valid on the rectangle, extends as a rational identity.  The second
factorization is treated identically.
With \(P_0=1/A\) and \(Q_0=1/C\), elementary algebra gives
\begin{equation}
 P_0(S-y)=R_0(S)+R_1(S)Q_0(y),
 \qquad R_1\ne0.
 \label{eq:translation-rank-prelim}
\end{equation}

We use the following elementary translation fact.  If rational
functions satisfy
\(P(S-y)=R_0(S)+R_1(S)Q(y)\), with \(R_1\ne0\), then either
\(P\) is constant or both \(P,Q\) are affine linear.  Indeed, all
translates \(P(S-\cdot)\) lie in the two-dimensional space
\(\operatorname{span}\{1,Q\}\).  A finite pole would produce
infinitely many translated pole locations, so \(P\) is a polynomial;
the translates of a degree-\(m\) polynomial span a space of dimension
\(m+1\).

If \(P_0\) is nonconstant, substitute affine expressions for
\(P_0,Q_0\) in \eqref{eq:translation-rank-prelim}.  Comparing
coefficients gives, after translating the variables,
\[
 r(x,y)=-\frac{y-y_0}{x-x_0}.
\]
Thus
\((x-x_0)f_x+(y-y_0)f_y=0\).  Applying Euler's identity to each
homogeneous part after translation forces every positive-degree part
of \(f\) to vanish, a contradiction.

Hence \(P_0\), and therefore \(A\), is constant.  We obtain
\(r-1=b(S)\).  In \((x,S)\)-coordinates, with \(G=G_X\), this is
\[
 G_x=b(S)G_S.
\]
If \(b=0\), then \(G\) is a polynomial in \(S\).  Otherwise choose
locally \(P_1'=1/b\).  The vector field
\(\partial_x-b(S)\partial_S\) annihilates both \(G\) and
\(x+P_1(S)\), so
\[
 G(x,S)=H(x+P_1(S)).
\]
Fixing \(S\) shows that \(H\) is a polynomial.  If \(\deg H=1\),
polynomiality of \(G=H(x+P_1(S))\) directly implies that \(P_1\) is a
polynomial.  If \(\deg H\ge2\), comparison of the next-to-leading
coefficient in \(x\) gives the same conclusion.  An affine rescaling gives
\(f(x,y)=H(ax+P(x+y))\), as required.
\end{proof}

\begin{corollary}[Three-coordinate classification]
\label{cor:three-coordinate}
The three polynomials \(f,G_X,G_Y\) are all analytic-special if and
only if \(f\in\mathfrak E_3\).  Consequently, if
\(f\notin\mathfrak E_3\), then at least one of the original or the
two sum-coordinate representations is not analytic-special.
\end{corollary}

\begin{proof}
The constant case is immediate, so assume that \(f\) is nonconstant.
By \cref{thm:two-transition}, write
\(f=H(ax+P(S))\).  The cases \(a=0\), \(P'\equiv0\), and
\(a+P'\equiv0\) immediately give a polynomial in one affine linear
form.  In the remaining cases,
\[
 \frac{f_x}{f_y}=\frac{a+P'(S)}{P'(S)}.
\]
If \(f\) is also analytic-special, this rational function of
\(x+y\) factors as \(A(x)B(y)\).  Its translates span a
one-dimensional space; the pole-and-translation argument used above
shows that it is constant.  Hence \(P'\) is constant and the argument
of \(H\) is affine linear.  The converse is immediate.
\end{proof}

\begin{corollary}[Directional-derivative test]
\label{cor:E3-test}
A nonconstant \(f\in\mathbb R[x,y]\) lies in \(\mathfrak E_3\) if and
only if there is a nonzero constant vector \((u,v)\) such that
\begin{equation}
 uf_x+vf_y\equiv0.
 \label{eq:E3-directional-test}
\end{equation}
\end{corollary}

\begin{proof}
If \eqref{eq:E3-directional-test} holds, take linear coordinates
parallel and perpendicular to \((u,v)\).  The polynomial is
independent of the parallel coordinate, so it is a polynomial in one
affine linear form.  The converse follows by differentiation.
\end{proof}

\begin{lemma}[A jointly Frostman shear estimate]
\label{lem:joint-shear}
Let \((X,Y)\) be jointly \((s,s,C)\)-Frostman with \(s>1/2\), and
put \(S=X+Y\).  For \(U=X\) or \(U=Y\),
\begin{equation}
 \mathbb P(U\in I,S\in J)
 \lesssim_C |I|^s|J|^{2s-1}
 \label{eq:joint-shear}
\end{equation}
for all intervals \(I,J\) of length at most one.
\end{lemma}

\begin{proof}
If \(|I|\le|J|\), the inverse shear of \(I\times J\) lies in a
bounded enlargement of a rectangle with side lengths \(|I|,|J|\).
Joint Frostman gives \(O(|I|^s|J|^s)\), which is stronger than
\eqref{eq:joint-shear}.  If \(|I|>|J|\), subdivide \(I\) into
\(O(|I|/|J|)\) intervals of length \(|J|\).  Each inverse image is
covered by a bounded number of \(|J|\)-squares, and hence has mass
\(O(|J|^{2s})\).  Summing gives
\(O(|I||J|^{2s-1})\le O(|I|^s|J|^{2s-1})\).
\end{proof}

The loss from \(s\) to \(2s-1\) in
\cref{lem:joint-shear} is genuine at the level of rectangular
Frostman information.  In particular, a direct transition argument
does not prove the jointly Frostman part of \Cref{thm:RZ-general} when
\(s\le1/2\).  The global input \Cref{prop:RZ-global} and the
localization developed below are used to avoid that loss.

\section{Entropy, algebraic, and localization preliminaries}
\label{sec:preliminaries}

\subsection{Dyadic entropy and Frostman consequences}

We first record elementary facts that will be used without further
comment.  For every fixed bounded \(Z\in\mathbb R^k\) and integer \(a\ge0\),
\begin{equation}
 H_n(Z)\le H_{n+a}(Z)\le H_n(Z)+ka.
 \label{eq:entropy-refinement}
\end{equation}
If \(T:\mathbb R^k\to\mathbb R^k\) is a fixed invertible linear map, then
\begin{equation}
 H_n(TZ)=H_n(Z)+O_T(1).
 \label{eq:linear-entropy}
\end{equation}
If \(F\) is fixed and Lipschitz on the range of \(Z\), then
\begin{equation}
 H_n(F(Z))\le H_n(Z)+O_F(1).
 \label{eq:Lipschitz-entropy}
\end{equation}
The same conclusions hold when the two variables differ by
\(O(2^{-n})\) almost surely.  They follow by observing that a fixed
enlargement of a level-\(n\) cell meets only \(O(1)\) level-\(n\)
cells.

\begin{lemma}[Bounded dyadic refinement]
\label{lem:bounded-refinement}
Let \(F_1,\ldots,F_m\) be fixed Lipschitz functions on a fixed
bounded subset of \(\mathbb R^k\).  The partition generated by
\[
 D_nZ,D_nF_1(Z),\ldots,D_nF_m(Z)
\]
has at most \(O(1)\) atoms over each atom of \(D_nZ\).  Its entropy
is at most \(H_n(Z)+O(1)\).
\end{lemma}

\begin{lemma}[Entropy supplied by Frostman bounds]
\label{lem:Frostman-entropy}
If \(X\) is \((s,C)\)-Frostman, then
\[
 H_n(X)\ge sn-\log C-O(1).
\]
If \((X,Y)\) is jointly \((s_1,s_2,C)\)-Frostman, then
\[
 H_n(X,Y)\ge(s_1+s_2)n-\log C-O(1).
\]
If the pair is conditionally
\((s_1,C_1;s_2,C_2)\)-Frostman, then also
\begin{align*}
 H(D_nX\mid D_nY)&\ge s_1n-O(1),\\
 H(D_nY\mid D_nX)&\ge s_2n-O(1).
\end{align*}
\end{lemma}

\begin{proof}
Every atom of the relevant dyadic partition has probability bounded
by a fixed constant times \(2^{-sn}\), or
\(2^{-(s_1+s_2)n}\) in the joint case.  Shannon entropy dominates
min-entropy.  For the conditional assertions, apply the conditional
atom bound on each conditioning cell and average.
\end{proof}

We also use the following elementary information inequalities.  If
\(L\) is a finite label, then
\begin{align}
 H(Z\mid L)&\le H(Z),
 \label{eq:label-entropy}\\
 \operatorname{I}(X;Y\mid L)&\le \operatorname{I}(X;Y)+H(L).
 \label{eq:label-MI}
\end{align}
For a mixture with weights \(w_j\), entropy concavity gives
\begin{equation}
 H(Z)\ge\sum_jw_jH(Z\mid L=j).
 \label{eq:mixture-concavity}
\end{equation}

\subsection{Polynomial sublevel sets}

The next lemma is the only place where the conditional Frostman
hypothesis is used geometrically.

\begin{lemma}[Polynomial Jacobian sublevel bound]
\label{lem:polynomial-sublevel}
Let \((X,Y)\) be bounded and conditionally
\((s_1,C_1;s_2,C_2)\)-Frostman, where \(s_1,s_2>0\).  If
\(R\in\mathbb R[x,y]\) is nonzero, then there are \(A<\infty\) and
\(\theta>0\) such that
\begin{equation}
 \mathbb P(|R(X,Y)|\le r)\le Ar^\theta,
 \qquad 0<r\le1.
 \label{eq:polynomial-sublevel}
\end{equation}
\end{lemma}

\begin{proof}
On a fixed compact neighborhood of the support, the polynomial
{\L}ojasiewicz inequality
{\cite[Corollary~2.6.7]{BCR}}
places \(\{|R|\le r\}\) in an
\(O(r^\gamma)\)-neighborhood of the real algebraic set \(Z(R)\),
for some \(\gamma>0\).  Semialgebraic stratification, followed by a
finite subdivision according to tangent direction, covers the
portion of \(Z(R)\) in that neighborhood by finitely many points and
finitely many Lipschitz graphs over one of the coordinate axes.
The \(t\)-neighborhood of a graph \(y=g(x)\) has conditional mass
\(O(C_2t^{s_2})\): condition on \(X=x\) and apply the conditional
Frostman bound to the interval centered at \(g(x)\).  This use of a
moving interval is justified by the disintegration observation after
the Frostman definition.  Graphs over the
other axis are treated with the other conditional bound, and points
are controlled by either marginal.  Taking the minimum of the
finitely many positive exponents proves \eqref{eq:polynomial-sublevel}.
\end{proof}

For a merely jointly \((s,s,C)\)-Frostman law one still has, when
\(s>1/2\),
\begin{equation}
 \mathbb P(\operatorname{dist}((X,Y),Z(R))\le r)\lesssim_R r^{2s-1}.
 \label{eq:joint-tube}
\end{equation}
Indeed, a Lipschitz arc may be covered by \(O(r^{-1})\) squares of
side \(r\), each of mass \(O(r^{2s})\).  This estimate explains why
a purely local proof for jointly Frostman laws changes character at
\(s=1/2\).

\subsection{Polynomial coordinate maps preserve entropy}

The following recovery statement is useful both for the direct
conditional argument and for entropy-profile preparation.  It is an
entropy comparison, not an entropy gain.

\begin{theorem}[Polynomial coordinate entropy]
\label{thm:coordinate-entropy}
Let \(L(x,y)=ax+by\) be nonzero, let \(P\in\mathbb R[x,y]\), and suppose
\begin{equation}
 J_{L,P}:=\det D(L,P)=aP_y-bP_x\not\equiv0.
 \label{eq:coordinate-J}
\end{equation}
If \((X,Y)\) is bounded and conditionally
\((s_1,C_1;s_2,C_2)\)-Frostman with \(s_1,s_2>0\), then
\begin{equation}
 \bigl|H_n(X,Y)-H_n(L(X,Y),P(X,Y))\bigr|\le B
 \label{eq:coordinate-recovery}
\end{equation}
for every \(n\ge1\), where \(B\) depends only on
\(L,P,s_1,s_2,C_1,C_2\) and the fixed support box, and in particular
is uniform over the admissible laws and independent of \(n\).  Hence
\begin{equation}
 \max\{H_n(L(X,Y)),H_n(P(X,Y))\}
 \ge\frac{s_1+s_2}{2}n-O(1).
 \label{eq:coordinate-baseline}
\end{equation}
For \(L=x+y\), condition \eqref{eq:coordinate-J} fails precisely
for \(P(x,y)=Q(x+y)\).
\end{theorem}

\begin{proof}
The Lipschitz inequality
\(H_n(L,P)\le H_n(X,Y)+O(1)\) follows from
\eqref{eq:Lipschitz-entropy}.  Complete \(L\) to a fixed invertible
linear coordinate system \((u,v)\).  In these coordinates the map is
\[
 F(u,v)=(u,\widetilde P(u,v)),
 \qquad \det DF=\partial_v\widetilde P,
\]
up to a fixed nonzero factor.

Let \(d=\deg_v\widetilde P\).  Except for finitely many values of
\(u\), the equation \(\widetilde P(u,v)=w\) has at most \(d\) real
solutions.  The exceptional values are common zeros of the
coefficients of the positive powers of \(v\) in \(\widetilde P\).
Each exceptional
source fibre is an affine line \(L(x,y)=u_0\).  Such a line has
probability zero under conditional Frostman regularity: use a
marginal bound for a vertical or horizontal line, and otherwise
condition on one original coordinate and shrink an interval in the
other coordinate.  The geometric
multiplicity of \(F\) is therefore at most \(d\) for almost every
target point.

Set \(K=0\) on \(\{|\partial_v\widetilde P|>1/2\}\), set \(K=k\)
on
\[
 2^{-k-1}<|\partial_v\widetilde P(U,V)|\le2^{-k},
 \qquad 1\le k<n,
\]
and set \(K=n\) when
\(|\partial_v\widetilde P(U,V)|\le2^{-n}\).  The polynomial
\(\partial_v\widetilde P(U,V)\) is a fixed nonzero multiple of
\(J_{L,P}(X,Y)\).  We apply \cref{lem:polynomial-sublevel} to the
latter polynomial in the original variables; no Frostman property
after the linear change \((X,Y)\mapsto(U,V)\) is being asserted.
Thus
\begin{equation}
 \mathbb P(K\ge k)\le A2^{-\theta k}\qquad(k\ge1).
 \label{eq:Jacobian-shell-tail}
\end{equation}
Consequently \(H(K)+\mathbb E K=O(1)\).

Fix a level-\(n\) output square and a shell \(k<n\).  If
\(k\le n-C_0\), then the Lipschitz variation of the polynomial
Jacobian on a fixed enlargement of any relevant input square is much
smaller than \(2^{-k}\), provided \(C_0\) is fixed sufficiently
large.  Thus \(|\det DF|\gtrsim2^{-k}\) on these enlarged squares.
Their images lie in a fixed enlargement of the selected output
square.  Bounded overlap, the area formula, and the a.e. multiplicity
bound give
\[
 \#\{\text{relevant input squares}\}\,2^{-2n-k}
 \lesssim \int|\det DF|
 \lesssim d\,2^{-2n}.
\]
There are therefore \(O(d2^k)\) possible input squares.  When
\(n-C_0<k<n\), the first output coordinate restricts \(u\) to
\(O(1)\) columns, leaving \(O(2^n)=O(2^k)\) input squares.  On the
terminal shell the bounded support gives the trivial \(O(2^{2n})\)
bound.

The chain rule now yields
\begin{align*}
 H(D_n(U,V)\mid D_nF(U,V))
 &\le H(K)+O(1)+\mathbb E[K\mathbf1_{K<n}]\\
 &\quad+2n\mathbb P(K=n)=O(1),
\end{align*}
where the terminal term is bounded by
\(O(n2^{-\theta n})\).  Returning through the fixed linear change
proves \eqref{eq:coordinate-recovery}.  Conditional Frostman implies
joint Frostman, so \cref{lem:Frostman-entropy} and subadditivity of
the two output coordinates prove \eqref{eq:coordinate-baseline}.
Finally, when \(L=x+y\), the equation \(P_y-P_x=0\) says in the
coordinates \(u=x+y,v=x-y\) that \(P_v=0\), hence
\(P=Q(u)\).
\end{proof}

\subsection{Simultaneous information profiles}

The nonlinear-projection inputs are statements about sets, whereas
the main theorems are statements about entropy.  The next lemma is
the conversion used throughout.  Its simultaneous form is important:
all source and image estimates hold on the same pieces.

\begin{lemma}[Simultaneous information profiles]
\label{lem:information-profiles}
Fix \(0<s<1\), numbers \(e,B\ge0\), a bounded box, and Lipschitz
scalar functions \(F_1,\ldots,F_m:\mathbb R^2\to\mathbb R\) on that box.  Let
\((X,Y)\) be jointly
\((s,s,C)\)-Frostman, put \(\delta=2^{-n}\), and suppose that
\begin{align}
 H_n(X,Y)&\le(2s+e)n+B,
 \label{eq:profile-source-upper}\\
 H_n(X),H_n(Y)&\le(s+e)n+B.
 \label{eq:profile-marginal-upper}
\end{align}
Fix \(\gamma,\zeta>0\) and \(K>1\).  Then there is a collection
\(\mathcal P_n\) of \(O_{m,K,\gamma,\zeta}(1)\) disjoint profiles
whose total probability \(W_n\) satisfies
\begin{equation}
 W_n\ge
 1-O\!\left(\frac{e}{\gamma}+\frac{B+1}{\gamma n}
              +\frac{m}{K}\right)-o_{n\to\infty}(1).
 \label{eq:profile-good-mass}
\end{equation}
Every profile \(\pi\in\mathcal P_n\) has probability at least
\(\delta^\zeta\).  Associated to it is a union \(E_\pi\) of
level-\(n\) source squares such that
\begin{align}
 \mathcal N_\delta(E_\pi)&\gtrsim\delta^{-2s+2\zeta},
 \label{eq:profile-set-size}\\
 \mathcal N_\delta(E_\pi\cap B(z,r))
 &\lesssim r^{2s}\delta^{-2s-\gamma-2\zeta}
 \qquad(\delta\le r\le1).
 \label{eq:profile-source-NC}
\end{align}
If \(A_\pi=x(E_\pi)\) and \(B_\pi=y(E_\pi)\), interpreted as
unions of level-\(n\) intervals, then
\begin{equation}
 |A_\pi\cap J|,\ |B_\pi\cap J|
 \lesssim |J|^s\delta^{1-s-\gamma-\zeta}
 \qquad(\delta\le |J|\le1).
 \label{eq:profile-marginal-NC}
\end{equation}
Finally, for each \(i\) there is a number \(a_i(\pi)\in[0,K]\)
such that
\begin{align}
 \mathcal N_\delta(F_i(E_\pi))
 &\lesssim\delta^{-a_i(\pi)-\zeta},
 \label{eq:profile-image-upper}\\
 H_n(F_i(X,Y))
 &\ge n\sum_{\pi\in\mathcal P_n}
          \mathbb P(\pi)a_i(\pi).
 \label{eq:profile-entropy-average}
\end{align}
All implicit constants may depend on the fixed data but not on
\(n\).
\end{lemma}

\begin{proof}
Refine the level-\(n\) source partition by the level-\(n\) bins of
all the functions \(F_i\).  By
\cref{lem:bounded-refinement}, each source square contains only
\(O(1)\) refined atoms and the refined source entropy is at most
\(H_n(X,Y)+O(1)\).

For a sampled refined atom, let \(L_0\) be minus the logarithm of
its probability.  Let \(L_x,L_y\) be the information of its global
\(X\)- and \(Y\)-bins, and let \(L_i\) be the information of its
global \(F_i\)-bin.  Joint and marginal Frostman bounds give
\begin{equation}
 L_0\ge2sn-O(1),\qquad L_x,L_y\ge sn-O(1).
 \label{eq:profile-information-lower}
\end{equation}
Their expectations are the corresponding entropies, up to \(O(1)\)
for \(L_0\).  Shifted Markov, using
\eqref{eq:profile-source-upper}--\eqref{eq:profile-marginal-upper},
therefore shows that restricting to
\begin{equation}
 L_0\le(2s+\gamma)n,\qquad
 L_x,L_y\le(s+\gamma)n
 \label{eq:profile-information-upper}
\end{equation}
discards at most
\(O(e/\gamma)+O((B+1)/(\gamma n))\) probability.  Since every
fixed bounded scalar output has entropy at most \(n+O(1)\), ordinary
Markov shows that restricting also to \(L_i\le Kn\) for all \(i\)
discards \(O(m/K)+O(m/(Kn))\) probability.

Divide every retained information range into half-open intervals of
length \(\zeta n\), and group atoms according to all of their
information intervals.  There are
\(O_{m,K,\gamma,\zeta}(1)\) resulting profiles.  Discard those of
probability below \(\delta^\zeta\); their total probability is
\(o(1)\).  This proves \eqref{eq:profile-good-mass}.

Fix a retained profile of mass \(w\ge\delta^\zeta\).  From every
source square occupied by the profile choose one of its refined
atoms, and let \(E_\pi\) be the union of these source squares.
Every chosen atom has mass \(O(\delta^{2s})\), and only \(O(1)\)
profile atoms can lie over a source square.  Hence
\(
 \mathcal N_\delta(E_\pi)\gtrsim w\delta^{-2s},
\)
which implies \eqref{eq:profile-set-size} for all sufficiently large
\(n\).

By \eqref{eq:profile-information-upper}, each chosen refined atom
has mass at least \(\delta^{2s+\gamma}\).  The joint Frostman bound
on a fixed enlargement of \(B(z,r)\) therefore gives
\eqref{eq:profile-source-NC}, with the displayed slack absorbing
cell-boundary constants.  Likewise every selected \(X\)-bin and
\(Y\)-bin has probability at least \(\delta^{s+\gamma}\).
Dividing the Frostman mass of an enlarged interval \(J\) by this
least bin mass, and then multiplying the resulting number of bins by
\(\delta\), proves \eqref{eq:profile-marginal-NC}.

If the \(L_i\)-window of the profile is
\([a_i n,(a_i+\zeta)n)\), every selected output bin has probability
at least \(\delta^{a_i+\zeta}\).  There are consequently at most
\(\delta^{-a_i-\zeta}\) such bins.  Lipschitzness adds only a fixed
number of neighboring bins when a selected source square is filled
out, proving \eqref{eq:profile-image-upper}.  Finally
\(L_i\ge a_i n\) on the profile; summing over the retained profiles
inside \(\mathbb E L_i=H_n(F_i(X,Y))\) gives
\eqref{eq:profile-entropy-average}.
\end{proof}

\begin{remark}[Finite refinements of profiles]
\label{rem:profile-refinement}
The construction in \cref{lem:information-profiles} is stable under
adjoining a label with a fixed finite range, even when the label
depends on \(n\).  Discard child profiles of mass below
\(\delta^{2\zeta}\).  Their total mass is \(o(1)\), and every
remaining child determines a union \(E\) of occupied source squares
satisfying \eqref{eq:profile-set-size}--\eqref{eq:profile-marginal-NC}
and the parent image bounds \eqref{eq:profile-image-upper}.
For the size estimate, each source square has mass
\(O(\delta^{2s})\), and the child has mass at least
\(\delta^{2\zeta}\).  For nonconcentration, choose in each occupied
square one refined atom from its parent profile.  These atoms are
disjoint across distinct squares, each has mass at least
\(\delta^{2s+\gamma}\), and those arising from a ball lie in a fixed
enlargement of that ball.  Dividing the Frostman bound by this atom
mass proves the required cell count.  The parent marginal and output
bin estimates pass to a subset.  Finally, the information lower
endpoints are inherited by each child, so their weighted averages
remain lower bounds for the original output entropies.  No Frostman
property of the child-conditioned law is required.
\end{remark}

\subsection{Buffered localization away from algebraic sets}

\begin{lemma}[Finite algebraic localization]
\label{lem:algebraic-localization}
Let \((X,Y)\) have bounded range and be conditionally
\((s_1,C_1;s_2,C_2)\)-Frostman, where \(s_1,s_2>0\).  Let
\(R_1,\ldots,R_\ell\in\mathbb R[x,y]\) be nonzero.  For every
\(\beta>0\) there is a finite rectangular grid label
\(L\in\{0,1,\ldots,M\}\) with the following properties:
\begin{enumerate}[label=\textup{(\roman*)}]
\item \(\mathbb P(L=0)\le\beta\).
\item For \(j\ge1\), the event \(\{L=j\}\) is a product rectangle
and is contained in the core of a buffered product rectangle on
which
\begin{equation}
 \min_{1\le\nu\le\ell}|R_\nu|\ge c_\beta>0.
 \label{eq:localization-lower-bound}
\end{equation}
\item If \(w_j=\mathbb P(L=j)>0\), the conditional law \(\mu_j\) on that
rectangle has marginal and rectangular estimates
\begin{align}
 \mu_j(X\in I)&\le C_1w_j^{-1}|I|^{s_1},\\
 \mu_j(Y\in J)&\le C_2w_j^{-1}|J|^{s_2},\\
 \mu_j(X\in I,Y\in J)
 &\le C_1C_2w_j^{-1}|I|^{s_1}|J|^{s_2}.
 \label{eq:localized-Frostman}
\end{align}
If a fixed smooth scalar map \(F\) has, on the buffered rectangle,
level sets uniformly transverse to the \(Y\)-axis, then also
\begin{equation}
 \mu_j(F\in I)\lesssim_F C_2w_j^{-1}|I|^{s_2};
 \label{eq:localized-output-cap}
\end{equation}
the symmetric assertion holds with \(X\) and \(Y\) interchanged.
\end{enumerate}
For a jointly \((s,s,C)\)-Frostman law with \(s>1/2\), conclusions
\textup{(i)}--\textup{(ii)} remain valid.  Moreover, on every good
rectangle of mass \(w_j>0\),
\begin{align}
 \mu_j(X\in I)&\lesssim Cw_j^{-1}|I|^s,\notag\\
 \mu_j(Y\in J)&\lesssim Cw_j^{-1}|J|^s,\notag\\
 \mu_j(X\in I,Y\in J)&\le Cw_j^{-1}|I|^s|J|^s.
 \label{eq:localized-Frostman-joint}
\end{align}
The implicit constants in the marginal estimates depend only on the
fixed support box.  In this case
\eqref{eq:localized-output-cap} is replaced by
\begin{equation}
 \mu_j(F\in I)\lesssim_F Cw_j^{-1}|I|^{2s-1}.
 \label{eq:localized-output-cap-joint}
\end{equation}
\end{lemma}

\begin{proof}
By \cref{lem:polynomial-sublevel}, a sufficiently small fixed
neighborhood of
\(Z(R_1\cdots R_\ell)\) has conditional-law probability at most
\(\beta/2\).  In the jointly Frostman case use
\eqref{eq:joint-tube}, which is effective when \(2s-1>0\).
Choose a fixed rectangular mesh much finer than the width of this
neighborhood.  Mark as bad every rectangle whose fixed enlargement
meets a smaller neighborhood of the algebraic set, and merge all bad
rectangles into the label \(0\).  With the two neighborhood widths
chosen in the proper order, their union lies in the original
neighborhood and therefore has probability at most \(\beta\).
Compactness supplies the uniform lower bound
\eqref{eq:localization-lower-bound} on each good buffered rectangle.

After conditioning on a good label of mass \(w_j\), divide the
original marginal or joint bound by \(w_j\) to obtain
\eqref{eq:localized-Frostman} or
\eqref{eq:localized-Frostman-joint}, respectively.  If the level sets of \(F\) are
transverse to the \(Y\)-axis, then, for each fixed \(x\), the
condition \(F(x,Y)\in I\) restricts \(Y\) to a uniformly bounded
number of intervals of length \(O_F(|I|)\).  Apply the conditional
Frostman estimate for \(Y\) given \(X=x\), justified by the
disintegration observation following the definition, integrate in
\(x\), and
divide by \(w_j\) to prove
\eqref{eq:localized-output-cap}.
For a jointly Frostman law, cover the \(r\)-neighborhood of each
uniformly regular level-set arc by \(O(r^{-1})\) squares of side
\(r=|I|\).  Dividing the resulting \(O(Cr^{2s-1})\) estimate by
\(w_j\) proves \eqref{eq:localized-output-cap-joint}.
\end{proof}

\begin{remark}[Averaging localized estimates]
\label{rem:localized-averaging}
The localized laws in \cref{lem:algebraic-localization} are used
through their marginal, rectangular, and transverse-output bounds;
no assertion that an arbitrary chart-conditioned law remains
conditionally Frostman is needed.  The losses are harmless under
averaging because
\begin{equation}
 \sum_{j=1}^M w_j\log\frac{1}{w_j}\le H(L)\le\log(M+1)
 \label{eq:chart-log-loss}
\end{equation}
and, for every discretized random variable \(Z\),
\begin{equation}
 H(Z)\ge\sum_{j=1}^M w_jH(Z\mid L=j).
 \label{eq:chart-entropy-average}
\end{equation}
For example, if every good chart satisfies
\(
 H_n(F\mid L=j)\ge an-B_0-B_1\log(1/w_j),
\)
then
\begin{equation}
 H_n(F)\ge(1-\beta)an-O_{B_0,B_1,M}(1).
 \label{eq:localized-globalization}
\end{equation}
\end{remark}

The final localization issue is different: a jointly Frostman law
can give substantial mass to a thin neighborhood of a non-axis
curve.  Near the critical values of an outer one-variable function,
one therefore uses the following entropy estimate rather than simply
discarding that neighborhood.

\begin{lemma}[Critical-layer entropy]
\label{lem:critical-layer}
Let \((X,Y)\) be jointly \((s,s,C)\)-Frostman in a fixed bounded
product chart, let \(p,q\) be fixed one-variable diffeomorphisms on
the coordinate intervals, and let \(T\subset\mathbb R\) be finite.  Put
\(\delta=2^{-n}\), \(R=\delta^\theta\) with \(0<\theta<1\), and
\begin{equation}
 \mathcal B_n=
 \{\operatorname{dist}(p(X)+q(Y),T)\le2R\},
 \qquad b_n=\mathbb P(\mathcal B_n).
 \label{eq:critical-layer}
\end{equation}
Then
\begin{equation}
 \frac{H_n(X)+H_n(Y)}{2}
 \ge sn+b_ns\theta n-h_2(b_n)-O(1),
 \label{eq:critical-layer-entropy}
\end{equation}
where \(h_2(t)=-t\log t-(1-t)\log(1-t)\).
\end{lemma}

\begin{proof}
For every level-\(n\) interval \(I\), bounded distortion of \(p,q\)
and rectangular Frostman regularity give
\begin{equation}
 \mathbb P(X\in I,\mathcal B_n)\lesssim\delta^sR^s,
 \qquad
 \mathbb P(Y\in I,\mathcal B_n)\lesssim\delta^sR^s.
 \label{eq:critical-layer-cell-cap}
\end{equation}
Indeed, after \(X\) is restricted to \(I\), the condition defining
\(\mathcal B_n\) restricts \(Y\) to a fixed finite union of
intervals of length \(O(R)\), and the other estimate is symmetric.
Thus
\begin{align*}
 H_n(X\mid\mathcal B_n)
 &\ge(s+s\theta)n+\log b_n-O(1),\\
 H_n(X\mid\mathcal B_n^c)
 &\ge sn+\log(1-b_n)-O(1).
\end{align*}
Entropy concavity over \(\mathcal B_n\) and its complement gives
\(
 H_n(X)\ge sn+b_ns\theta n-h_2(b_n)-O(1).
\)
The same argument for \(Y\), followed by averaging, proves
\eqref{eq:critical-layer-entropy}.  The cases \(b_n=0,1\) follow by
the usual convention \(0\log0=0\).
\end{proof}
\section{Proof of Theorem~\ref{thm:RZ-general}}
\label{sec:RZ-proof}

We prove the conditionally Frostman and jointly Frostman assertions
simultaneously.  The expanding-polynomial and flat four-web inputs
are Proposition~\ref{prop:RZ-global} and
Proposition~\ref{prop:RZ-flat}, respectively.

\subsection{Entropy profiles}

We specialize \cref{lem:information-profiles} to the coordinate
maps \(x,y\) and the observables \(x+y,f\), keeping the information-window
notation needed by the two nonlinear-projection inputs.  We record
the specialization explicitly because the later chart and near--far
labels must be attached to the same profiles.  This profile reduction is used in both parts of
Theorem~\ref{thm:RZ-general}.  Apply the same affine change to both
source coordinates so that the support is contained in
\([1/4,3/4]^2\subset[0,1]^2\).
This changes all dyadic entropies by \(O_c(1)\), replaces \(X+Y\) by
an affine function of the normalized sum, and replaces \(f\) by a
fixed affine precomposition.  We continue to write \(X,Y,f\) for the
normalized objects.  The dependence of \(\varepsilon\) on \(c\) allows
for this normalization.

Put
\[
 S=X+Y,\qquad \delta=2^{-n}.
\]
Since \(f\notin\mathfrak E_3\), one has
\(f_y-f_x\not\equiv0\); indeed, \(f_y-f_x\equiv0\) would imply
\(f(x,y)=Q(x+y)\).  The polynomial coordinate-recovery theorem,
Theorem~\ref{thm:coordinate-entropy}, applied to
\((x+y,f(x,y))\), therefore
gives
\begin{equation}
 H_n(X,Y)=H_n(S,f(X,Y))+O_{f,s,C,c}(1)
 \label{eq:RZ-entropy-recovery}
\end{equation}
under the conditional Frostman hypothesis.

Suppose first that \eqref{eq:RZ-general-cond} fails at
level \(n\) with a positive number \(\varepsilon\).  The conditional
Frostman entropy bounds and \eqref{eq:RZ-entropy-recovery} imply
\begin{align}
 H_n(X,Y)&\le2(s+\varepsilon)n+O(1),
 \label{eq:RZ-source-upper}\\
 H_n(X),H_n(Y)&\le(s+2\varepsilon)n+O(1).
 \label{eq:RZ-marginal-upper}
\end{align}
For example,
\(H(D_nY\mid D_nX)\ge sn-O(1)\), and hence
\(H_n(X)\le H_n(X,Y)-sn+O(1)\).

If \eqref{eq:RZ-general-joint} fails, then
\[
 H_n(X)+H_n(Y)\le2(s+\varepsilon)n.
\]
The marginal Frostman lower bounds give
\eqref{eq:RZ-marginal-upper}, while subadditivity gives
\eqref{eq:RZ-source-upper}.  In addition,
\begin{equation}
 H_n(X,Y)\ge2sn-O(1).
 \label{eq:RZ-source-lower}
\end{equation}
Thus \eqref{eq:RZ-source-upper}--\eqref{eq:RZ-marginal-upper}
hold under either failed conclusion.  Under conditional Frostman,
\eqref{eq:RZ-source-lower} also holds because conditional Frostman
implies joint Frostman.

Apply \cref{lem:information-profiles} with \(e=2\varepsilon\), with \(B\)
absorbing the fixed errors above, and with the four scalar maps
\[
 \pi_x=x,\qquad \pi_y=y,\qquad \pi_S=x+y,\qquad \pi_f=f.
\]
Choose parameters \(0<\gamma,\zeta<1\) and \(K>1\).
The source and marginal information truncation discards mass at most
\begin{equation}
 O\!\left(\frac{\varepsilon}{\gamma}\right)
 +O\!\left(\frac{1}{\gamma n}\right),
 \label{eq:RZ-profile-tail-loss}
\end{equation}
and the output truncation and small-profile deletion cost
\(O(K^{-1})+o(1)\).  Thus a collection of a fixed number of profiles
retains mass \(1-O(\varepsilon/\gamma+K^{-1})-o(1)\).

For each retained profile, write \(a_i n\) for the lower endpoint
of its global \(\pi_i\)-information window, where
\(i\in\{x,y,S,f\}\), and let \(E\) be its union of occupied
level-\(n\) source squares.  Put \(A=x(E)\) and \(B=y(E)\).
The conclusions of the profile lemma, expressed as area and length
estimates, are
\begin{equation}
 |E|\gtrsim\delta^{2-2s+2\zeta},
 \label{eq:RZ-profile-size}
\end{equation}
\begin{equation}
 |A\cap J|,|B\cap J|
 \lesssim |J|^s\delta^{1-s-\gamma-\zeta}
 \qquad(\delta\le |J|\le1),
 \label{eq:RZ-profile-marginal-nc}
\end{equation}
and
\begin{equation}
 |E\cap B(z,r)|
 \lesssim r^{2s}\delta^{2-2s-\gamma-2\zeta}
 \qquad(\delta\le r\le1).
 \label{eq:RZ-profile-source-nc}
\end{equation}
For \(|J|<\delta\), the marginal bound needed by the global input
follows from \(|A\cap J|,|B\cap J|\le |J|\); for \(|J|>1\),
it follows from the same bound at length one and the fixed support.
The image bounds are
\begin{equation}
 |f(E)|\lesssim\delta^{1-a_f-\zeta},\qquad
 \mathcal N_\delta(\pi_i(E))\lesssim\delta^{-a_i-\zeta}
 \quad(i=x,y,S,f),
 \label{eq:RZ-profile-image-upper}
\end{equation}
where the length bound follows by covering the scalar image by
intervals.  The weighted averages of the \(a_i n\) are lower bounds
for the corresponding entropies by
\eqref{eq:profile-entropy-average}.  When fixed finite chart or
near--far labels are adjoined below, we use
\cref{rem:profile-refinement}, retaining children of mass at least
\(\delta^{2\zeta}\).  All estimates above and the information
lower endpoints then remain available on the same pieces.

\subsection{The non-special case}

Assume that \(f\) is not a polynomial special form.  Let
\(\eta_{\rm RZ}(f,s)\) and \(\kappa_{\rm RZ}(s)\) be as in
Proposition~\ref{prop:RZ-global}.  Choose
\begin{equation}
 2\zeta<\eta_{\rm RZ},\qquad
 \gamma+\zeta<\eta_{\rm RZ},\qquad
 \zeta<\frac{1}{100}\kappa_{\rm RZ}.
 \label{eq:RZ-global-parameters}
\end{equation}
For all sufficiently large \(n\), fixed implicit constants are
absorbed by the strict inequalities in
\eqref{eq:RZ-global-parameters}.  Thus
\eqref{eq:RZ-profile-size} and
\eqref{eq:RZ-profile-marginal-nc} satisfy the hypotheses of
Proposition~\ref{prop:RZ-global}, and
\[
 |f(E)|\ge\delta^{1-s-\kappa_{\rm RZ}}.
\]
Comparing this with \eqref{eq:RZ-profile-image-upper} yields
\begin{equation}
 a_f\ge s+\kappa_{\rm RZ}-O(\zeta)
 \label{eq:RZ-profile-output-gain}
\end{equation}
for every retained profile.

Let \(W\) be the total mass of the retained profiles.  By
\eqref{eq:RZ-profile-tail-loss},
\begin{equation}
 W\ge1-O\!\left(\frac{\varepsilon}{\gamma}\right)
       -O_s(K^{-1})-o(1).
 \label{eq:RZ-retained-mass-global}
\end{equation}
Averaging the information lower endpoints in
\eqref{eq:RZ-profile-output-gain} gives
\[
 H_n(f(X,Y))
 \ge W\bigl(s+\kappa_{\rm RZ}-O(\zeta)\bigr)n-O(1).
\]
Fix \(\gamma,\zeta\) satisfying
\eqref{eq:RZ-global-parameters}; then take \(K\) sufficiently large,
then \(\varepsilon>0\) sufficiently small in terms of
\(\gamma\kappa_{\rm RZ}\), and finally \(n\) sufficiently large.
The last display contradicts
\(H_n(f(X,Y))\le(s+\varepsilon)n\).  This proves both parts of
Theorem~\ref{thm:RZ-general} when \(f\) is not a polynomial special
form.

\subsection{Polynomial special forms}

It remains to consider a polynomial special form
\(f\notin\mathfrak E_3\).  We isolate the chart construction because
the outer one-variable function is allowed to have critical points.

\begin{lemma}[Regular charts for polynomial special forms]
\label{lem:RZ-special-charts}
Fix a bounded source box and let \(f\) be a nonconstant polynomial
special form outside \(\mathfrak E_3\).  Suppose that the two source
marginals are \(s\)-Frostman.  For every \(\beta>0\), outside a union
of finitely many coordinate strips of total probability at most
\(\beta\), the source box is covered by finitely many buffered product
rectangles with the following property.  On each rectangle there are
one-variable real-analytic functions \(h,p,q,F,G\) such that
\begin{equation}
 f(x,y)=h(p(x)+q(y)),\qquad
 x+y=F(p(x))+G(q(y)),
 \label{eq:RZ-flat-normal-form}
\end{equation}
where \(p,q\) and their inverses \(F,G\) are local diffeomorphisms.
After interchanging the two coordinates if necessary, one also has
\begin{equation}
 F',F'',G'\ne0
 \label{eq:RZ-chart-regularity}
\end{equation}
throughout the buffered rectangle.  No nonvanishing assumption is
made on \(h'\).
\end{lemma}

\begin{proof}
Write first \(f(x,y)=H(A(x)+B(y))\).  Since
\(f\notin\mathfrak E_3\), both \(A\) and \(B\) are nonconstant.
Away from the finitely many zeros of \(A'\) and \(B'\), set
\[
 p=A,\qquad q=B,\qquad h=H,
 \qquad F=A^{-1},\qquad G=B^{-1}.
\]
This gives \eqref{eq:RZ-flat-normal-form}; for example,
\[
 F''=-\frac{A''\circ F}{(A'\circ F)^3}.
\]
If both \(F''\) and \(G''\) vanished identically on an open product
chart, then \(F\) and \(G\) would be affine there.  It would follow
that \(p(x)+q(y)\) is affine in \((x,y)\) on an open set.  A fixed
nonzero directional derivative of the polynomial \(f\) would then
vanish on that open set, hence identically, placing \(f\) in
\(\mathfrak E_3\), a contradiction.

Next write \(f(x,y)=H(A(x)B(y))\).  Besides the zeros of \(A'\) and
\(B'\), remove the zeros of \(A\) and \(B\).  On a resulting sign
rectangle let \(\sigma=\operatorname{sgn}(A(x)B(y))\) and set
\[
 p(x)=\ln|A(x)|,\qquad q(y)=\ln|B(y)|,
 \qquad h(t)=H(\sigma e^t),
\]
with \(F=p^{-1}\) and \(G=q^{-1}\).  Again
\eqref{eq:RZ-flat-normal-form} holds.  In this case \(F''\) and
\(G''\) cannot both vanish identically: otherwise \(p\) and \(q\)
would be affine on intervals, so the nonconstant polynomials \(A\)
and \(B\) would agree there with nonconstant exponential functions,
which is impossible.

Thus, in either case, at least one of \(F''\), \(G''\) is a
nonzero analytic function on every regular chart.  Its zeros on a
buffered compact subinterval are finite.  Delete small fixed
neighborhoods of those zeros, interchange the two coordinates when
the chosen function is \(G''\), and subdivide into finitely many
buffered product rectangles.  All deleted sets are finite unions of
vertical or horizontal strips.  Their total probability can be made
at most \(\beta\) by the marginal Frostman estimates.  The remaining
compact rectangles satisfy \eqref{eq:RZ-chart-regularity} uniformly.
\end{proof}

Apply \cref{lem:RZ-special-charts}.  Subdivide its rectangles further
so that, in \((p,q)\)-coordinates, every chart is a buffered product
of intervals of the same fixed length.  A common affine rescaling of
\(p\) and \(q\) then preserves \(p+q\), up to an affine output
change.  By \cref{lem:buffered-coordinate-change}, the size, ball
nonconcentration, and covering estimates survive the passage from
\((x,y)\) to \((p,q)\).  Proposition~\ref{prop:RZ-flat} therefore
applies on each chart to the four maps
\begin{equation}
 p,\qquad q,\qquad p+q,\qquad F(p)+G(q)=S.
 \label{eq:RZ-four-chart-maps}
\end{equation}

Adjoin the finite chart label to the profiles constructed above.
Whenever this label or the near--far label below splits a profile,
retain only children of mass at least \(\delta^{2\zeta}\).
By \cref{rem:profile-refinement}, the total additional loss is
\(o(1)\), and
\eqref{eq:RZ-profile-size}--\eqref{eq:RZ-profile-image-upper}
remain valid, with the parent's information lower endpoints.

We next control the critical points of the outer function \(h\).
The estimate below is the finite chart-labelled version of
\cref{lem:critical-layer}.
Since \(f\) is nonconstant, \(h'\) is not identically zero on any
chart.  On the finitely many buffered charts, let \(T\) be the
corresponding finite union of zero sets of \(h'\).  Fix
\(0<\theta<1\), put
\(R=\delta^\theta\), and define the chartwise event
\begin{equation}
 \mathcal B_n
 =\{\operatorname{dist}(p(X)+q(Y),T)\le2R\}.
 \label{eq:RZ-critical-layer}
\end{equation}
The notation in \eqref{eq:RZ-critical-layer} is interpreted on each
chart, and \(\mathcal B_n\) is the union of the finitely many
chart-labelled events.  Fixed derivative bounds and joint Frostman
give, for every level-\(n\) \(X\)-bin \(I\),
\begin{equation}
 \mathbb P(X\in I,\mathcal B_n)\lesssim\delta^sR^s,
 \label{eq:RZ-critical-X-bin}
\end{equation}
and the analogous estimate holds for \(Y\).  Indeed, after \(X\) is
restricted to \(I\), the condition in
\eqref{eq:RZ-critical-layer} restricts \(Y\) to finitely many
intervals of length \(O(R)\).  Conditional Frostman implies the same
joint Frostman bound, so \eqref{eq:RZ-critical-X-bin} holds in both
parts of the theorem.

Let \(q_n=\mathbb P(\mathcal B_n)\).  The atom caps in
\eqref{eq:RZ-critical-X-bin} and the ordinary marginal atom caps give
\begin{align*}
 H_n(X\mid\mathcal B_n)&\ge(s+s\theta)n+\log q_n-O(1),\\
 H_n(X\mid\mathcal B_n^c)&\ge sn+\log(1-q_n)-O(1),
\end{align*}
with the evident interpretation when \(q_n=0\) or \(1\).  Entropy
concavity, followed by the analogous estimate for \(Y\), gives
\begin{equation}
 \frac{H_n(X)+H_n(Y)}{2}
 \ge sn+q_ns\theta n-h_2(q_n)-O(1),
 \label{eq:RZ-critical-entropy}
\end{equation}
where \(h_2\) is the binary entropy.  Under failure of the joint
conclusion, \eqref{eq:RZ-critical-entropy} implies
\begin{equation}
 q_n=O\!\left(\frac{\varepsilon}{s\theta}\right)+o(1).
 \label{eq:RZ-critical-mass-joint}
\end{equation}
Under failure of the conditional conclusion,
\eqref{eq:RZ-marginal-upper} gives instead
\begin{equation}
 q_n=O\!\left(\frac{2\varepsilon}{s\theta}\right)+o(1).
 \label{eq:RZ-critical-mass-conditional}
\end{equation}

Refine the profiles by the event \(\mathcal B_n\), retain only the
far child profiles of mass at least \(\delta^{2\zeta}\), and saturate
their source atoms to level-\(n\) squares.  Since \(\delta=o(R)\), for
large \(n\) every such saturated square remains at distance at least
\(R\) from \(T\).  Fix one of the polynomial-special representations
used in \cref{lem:RZ-special-charts}, with outer polynomial \(H\), and
put
\[
 m=\max\{0,\deg H-1\}.
\]
This number depends only on the fixed polynomial representation, not
on the chart decomposition or on the Frostman constant.  In the
additive charts \(h'=H'\), while in the multiplicative charts
\(h'(t)=\sigma e^tH'(\sigma e^t)\).  Thus every chartwise zero of
\(h'\) has order at most \(m\).  On the far region,
\[
 |h'|\gtrsim R^m=\delta^{m\theta}.
\]
Splitting into the finitely many monotonicity intervals of \(h\)
therefore gives
\begin{equation}
 \mathcal N_\delta((p+q)(E))
 \lesssim\delta^{-m\theta}\mathcal N_\delta(f(E)).
 \label{eq:RZ-outer-critical-loss}
\end{equation}

For one retained far profile define
\begin{equation}
 b_i=\frac{\log \mathcal N_\delta(\pi_i(E))}{\log(1/\delta)},
 \qquad i\in\{x,y,S,f\}.
 \label{eq:RZ-covering-exponents}
\end{equation}
The maps \((x,y),(x,S),(y,S),(x,p+q),(y,p+q)\) are uniformly
bi-Lipschitz on each buffered chart.  Combining this with
\eqref{eq:RZ-profile-size} and
\eqref{eq:RZ-outer-critical-loss} gives
\begin{align}
 b_x+b_y,\ b_x+b_S,\ b_y+b_S
 &\ge2s-O(\zeta),
 \label{eq:RZ-pair-bounds-1}\\
 b_x+b_f,\ b_y+b_f
 &\ge2s-m\theta-O(\zeta).
 \label{eq:RZ-pair-bounds-2}
\end{align}

Use the uniform exponents in Proposition~\ref{prop:RZ-flat} for all
of the finitely many charts.  Thus
\(\eta_{\rm web},\kappa_{\rm web}>0\) depend only on \(s\); the
threshold scale and implicit constants may depend on the charts.
Choose \(\gamma,\zeta\) so that
\begin{equation}
 2\zeta<\eta_{\rm web},\qquad
 \gamma+2\zeta<\eta_{\rm web}.
 \label{eq:RZ-web-parameters}
\end{equation}
The size and nonconcentration estimates
\eqref{eq:RZ-profile-size} and
\eqref{eq:RZ-profile-source-nc}, transferred through the fixed chart
coordinates, meet the hypotheses of
Proposition~\ref{prop:RZ-flat}.  Put
\(\kappa_0=\kappa_{\rm web}/2\).  After taking \(n\)
sufficiently large, Proposition~\ref{prop:RZ-flat},
\eqref{eq:RZ-four-chart-maps}, and
\eqref{eq:RZ-outer-critical-loss} imply
\begin{equation}
 \max\{b_x,b_y,b_S,b_f\}
 \ge s+\kappa_0-m\theta.
 \label{eq:RZ-web-max}
\end{equation}
Here we also used the uniform one-dimensional comparisons
\[
 \mathcal N_\delta(p(E))\asymp\mathcal N_\delta(x(E)),
 \qquad
 \mathcal N_\delta(q(E))\asymp\mathcal N_\delta(y(E)),
\]
which follow from the bi-Lipschitz inverse functions \(F,G\).  If the
large image supplied by the flat theorem is \((p+q)(E)\),
\eqref{eq:RZ-outer-critical-loss} accounts for the loss
\(m\theta\) in \eqref{eq:RZ-web-max}.

We use the following elementary consequence of the pair bounds
\eqref{eq:RZ-pair-bounds-1}--\eqref{eq:RZ-pair-bounds-2} and the
maximum bound \eqref{eq:RZ-web-max}:
\begin{equation}
 \frac{3}{2}(b_x+b_y)+b_S+b_f
 \ge5s+\kappa_0
      -O(m\theta+\gamma+\zeta).
 \label{eq:RZ-profile-LP}
\end{equation}
For completeness, subtract \(s\) from the four variables and write
the resulting variables as \(x,y,z,w\).  With
\(e=O(\gamma+\zeta)\), the pair bounds become
\[
 x+y,\ x+z,\ y+z\ge-e,
 \qquad x+w,\ y+w\ge-m\theta-e.
\]
If \(x\) is the maximum, then
\[
 \frac{3}{2}(x+y)+z+w
 \ge\frac{3}{2}x-\frac{1}{2}y-m\theta-O(e)
 \ge x-m\theta-O(e),
\]
and the case in which \(y\) is the maximum is symmetric.  If \(z\)
is the maximum, put \(T_0=x+y\ge-e\).  Since
\(\min(x,y)\le T_0/2\), the same expression is at least
\(z+T_0-m\theta-O(e)\).  If \(w\) is the maximum, it is at least
\(w+T_0-O(e)\).  Combining these alternatives with
\eqref{eq:RZ-web-max} proves \eqref{eq:RZ-profile-LP}.

If \(a_x,a_y,a_S,a_f\) are the lower endpoints of the information
windows of the profile, \eqref{eq:RZ-profile-image-upper} gives
\begin{equation}
 b_i\le a_i+\zeta+O(1/n),
 \qquad i\in\{x,y,S,f\}.
 \label{eq:RZ-cover-information-bridge}
\end{equation}
All coefficients in \eqref{eq:RZ-profile-LP} are positive.  Hence
\eqref{eq:RZ-profile-LP} remains valid with \(a_i\) in place of
\(b_i\), after enlarging the \(O(\zeta)\) error.

Let \(W_{\mathrm{far}}\) be the total mass of the retained far
profiles.  Equations \eqref{eq:RZ-profile-tail-loss},
\eqref{eq:RZ-critical-mass-joint}, and
\eqref{eq:RZ-critical-mass-conditional} give, in either case,
\begin{equation}
 W_{\mathrm{far}}
 \ge1-O\!\left(\beta+\frac{\varepsilon}{\gamma}
 +\frac{\varepsilon}{s\theta}+\frac{1}{K}\right)-o(1).
 \label{eq:RZ-far-mass}
\end{equation}
Averaging the information form of \eqref{eq:RZ-profile-LP} over
these profiles yields
\begin{align}
 &\frac{3}{2}\bigl(H_n(X)+H_n(Y)\bigr)
   +H_n(S)+H_n(f(X,Y))
 \nonumber\\
 &\qquad\ge
 W_{\mathrm{far}}
 \bigl(5s+\kappa_0
       -O(m\theta+\gamma+\zeta)\bigr)n-O(1).
 \label{eq:RZ-averaged-LP}
\end{align}

Under failure of the conditional conclusion,
\eqref{eq:RZ-marginal-upper} gives the upper bound
\begin{equation}
 \frac{3}{2}\bigl(H_n(X)+H_n(Y)\bigr)
 +H_n(S)+H_n(f(X,Y))
 \le(5s+8\varepsilon)n+O(1).
 \label{eq:RZ-conditional-LP-upper}
\end{equation}
Under failure of the joint conclusion, the corresponding upper bound
is
\begin{equation}
 \frac{3}{2}\bigl(H_n(X)+H_n(Y)\bigr)
 +H_n(S)+H_n(f(X,Y))
 \le(5s+5\varepsilon)n.
 \label{eq:RZ-joint-LP-upper}
\end{equation}

Finally, choose the parameters in the following order.  The integer
\(m\) above is already fixed by the chosen polynomial-special
representation.  First take
\[
 \theta<\frac{\kappa_0}{100(m+1)}.
\]
Next choose \(\beta>0\) sufficiently small in terms of the uniform
web gain \(\kappa_0\), construct the finitely many buffered charts in
\cref{lem:RZ-special-charts}, and then choose
\(\gamma,\zeta,K^{-1}\) sufficiently small, with
\eqref{eq:RZ-web-parameters} in force, so that the errors and the
fixed mass loss in \eqref{eq:RZ-averaged-LP} use less than one half
of the web gain.  Next choose
\[
 \varepsilon\ll\min\{\gamma,s\theta,\kappa_0\},
\]
and finally take \(n\) sufficiently large.  The lower bound
\eqref{eq:RZ-averaged-LP} then contradicts
\eqref{eq:RZ-conditional-LP-upper} or
\eqref{eq:RZ-joint-LP-upper}, respectively.  This completes the proof
of Theorem~\ref{thm:RZ-general}. \hfill\(\square\)
\section{Explicit entropy gains from Pham's folding estimate}
\label{sec:Pham-explicit}

This section proves \cref{thm:Pham-explicit}.  The analytic input is the
two-sided folding estimate of Pham \cite{PhamFIO}.  The fold costs
\(1/6\) of a derivative and, together with the cross partition,
produces the threshold loss \(2/3\).  After localization and an
information-theoretic comparison with the product of the marginals,
this gives the explicit
entropy exponents below.  The algebraic choice of coordinates has
already been made in \cref{thm:two-transition,cor:three-coordinate}.

All implicit constants in this section are independent of the dyadic
scale \(n\) and of the admissible law.  They may depend on the fixed
polynomial, the support bound, the Frostman constants, the displayed
exponents, and any strict entropy exponent chosen in the proof.  As before, \(\log\) denotes the base-two logarithm.

\subsection{Entropy and energy lemmas}

For later use, note that \eqref{eq:entropy-refinement} and data
processing give, for every fixed integer \(r\geq0\),
\begin{equation}
 \operatorname{I}_n(U;V)\leq \operatorname{I}_{n+r}(U;V)
 \leq \operatorname{I}_n(U;V)+2r.
 \label{eq:Pham-information-refinement}
\end{equation}
Indeed, the first inequality is data processing.  For the second,
the chain rule bounds the extra information carried by each of
\(D_{n+r}U\) and \(D_{n+r}V\), beyond its level-\(n\) parent, by
\(r\).  We shall also use the fixed linear-change estimate
\eqref{eq:linear-entropy}, entropy concavity
\eqref{eq:mixture-concavity}, and \cref{lem:Frostman-entropy}.  In
particular, a jointly \((a,b,K)\)-Frostman pair satisfies
\begin{equation}
 H_n(U,V)\geq(a+b)n-O_K(1).
 \label{eq:Pham-joint-entropy-lower}
\end{equation}

For a compactly supported probability measure \(\mu\) on \(\mathbb R\)
and $0<a<1$, put
\[
 I_a(\mu)=\iint |x-y|^{-a}\,d\mu(x)d\mu(y).
\]
The layer-cake formula shows that if $\mu$ is $(s,C)$-Frostman and
$0<a<s$, then
\begin{equation}
 I_a(\mu)\ll_{a,s,c}1+C.
 \label{eq_pham_input_energy}
\end{equation}
If $\nu$ is a compactly supported probability measure and
$I_u(\nu)<\infty$, then
\begin{equation}
 H_n(\nu)\geq u n-\log I_u(\nu),\qquad 0<u<1.
 \label{eq_pham_energy_entropy}
\end{equation}
Indeed, on a dyadic interval $Q$ of length $2^{-n}$ one has
$|x-y|^{-u}\geq2^{un}$ for $x,y\in Q$.  Hence
\[
 I_u(\nu)\geq2^{un}\sum_Q\nu(Q)^2,
\]
and Shannon entropy is at least collision entropy.

We shall also use the following elementary inequality.  For a
probability vector $p=(p_j)$, write
\[
 H_{\mathrm{coll}}(p)=-\log\sum_jp_j^2,
 \qquad
 D(p\Vert q)=\sum_jp_j\log\frac{p_j}{q_j}.
\]

\begin{lemma}[Entropy and relative entropy]
\label{lem_pham_relative_entropy}
For any two probability vectors $p,q$,
\begin{equation}
 H(p)+2D(p\Vert q)\geq H_{\mathrm{coll}}(q).
 \label{eq_pham_relative_entropy}
\end{equation}
\end{lemma}

\begin{proof}
The assertion is immediate when $D(p\Vert q)=+\infty$.  Otherwise,
put $Z=\sum_jq_j^2$ and $r_j=q_j^2/Z$.  Then
\[
 H(p)+2D(p\Vert q)
 =\sum_jp_j\log\frac{p_j}{q_j^2}
 =D(p\Vert r)-\log Z\geq-\log Z.
\]
\end{proof}

\subsection{The folding input}

Let $f\in\mathbb R[x,y]$.  On the set where $f_{xy}\neq0$, define
\begin{equation}
 \kappa_f=\nabla f\wedge
 \nabla\left(\frac{f_xf_y}{f_{xy}}\right).
 \label{eq_pham_kappa}
\end{equation}
It is convenient to clear the denominator and put
\begin{align}
 \mathcal B_f
 &=f_y^2(f_xf_{xxy}-f_{xx}f_{xy})
 -f_x^2(f_yf_{xyy}-f_{xy}f_{yy}),
 \label{eq_pham_cleared_curvature}\\
 P_f&=f_xf_yf_{xy}\mathcal B_f.
 \label{eq_pham_bad_polynomial}
\end{align}
Here and below, $Z(P_f)=\{(x,y):P_f(x,y)=0\}$.
A direct calculation gives
\begin{equation}
 \mathcal B_f=f_{xy}^{2}\kappa_f.
 \label{eq_pham_curvature_identity}
\end{equation}
In this section, ``non-special'' means ``not analytic-special.''
We use analytic-specialness in the sense of
\eqref{eq:analytic-special}.  Multiplicative polynomial forms are
included locally after taking logarithms on a regular sign component.
If none of \(f_x,f_y,f_{xy}\) vanishes identically, then
\(\kappa_f\equiv0\) if and only if \(f\) is analytic-special; this is
\cite[Lemma~4.4]{PhamFIO}, which in turn cites
\cite[Lemma~3.3]{RZ}.  A non-special polynomial therefore has
\(f_xf_yf_{xy}\mathcal B_f\not\equiv0\), and hence \(P_f\) is a
nonzero polynomial.
\begin{proposition}[Local folding energy estimate]
\label{prop_pham_local_energy}
Let $f\in\mathbb R[x,y]$ be non-special.  Let
\[
 K=I\times J\Subset K^+=I^+\times J^+
 \Subset\mathbb R^2\setminus Z(P_f),
\]
where $I,J,I^+,J^+$ are compact intervals.  Suppose that $\mu_1$ and
$\mu_2$ are probability measures supported on $I$ and $J$,
respectively.  If
\[
 0<a_1,a_2,u<1,
 \qquad a_1+a_2>u+\frac{2}{3},
\]
then
\begin{equation}
 I_u\big(f_*(\mu_1\times\mu_2)\big)
 \leq A
 \bigl(1+I_{a_1}(\mu_1)\bigr)
 \bigl(1+I_{a_2}(\mu_2)\bigr),
 \label{eq_pham_local_energy}
\end{equation}
where $A$ depends on the displayed fixed data, but not on the two
measures.
\end{proposition}

\begin{proof}
By \eqref{eq_pham_curvature_identity}, the functions
$f_x,f_y,f_{xy},\kappa_f$ are bounded away from zero on $K^+$.
Put \(\rho=f_xf_y/f_{xy}\).  We first record why the regular
product rectangle is a valid chart for the fold argument.
The derivatives \(f_x,f_y\) have fixed nonzero signs on \(K^+\).
For every real \(t\), the level set \(f=t\) in the interior of
\(K^+\), when nonempty, is a connected graph \(y=y_t(x)\) over an
interval.  To see connectedness, assume for example that \(f_y>0\)
and write \(y_-,y_+\) for the vertical endpoints.  The permissible
\(x\)'s satisfy \(f(x,y_-)<t<f(x,y_+)\); both endpoint functions
are strictly monotone in the same direction, so this set is an
interval.  The other sign choices give the same conclusion.
Along this graph,
\begin{equation}
 \frac{d}{dx}\rho(x,y_t(x))
 =\rho_x-\frac{f_x}{f_y}\rho_y
 =-\frac{\kappa_f}{f_y}.
 \label{eq:Pham-level-monotonicity}
\end{equation}
This derivative has a fixed nonzero sign.  Consequently
\((f,\rho)\) is injective on the rectangle's interior; its Jacobian
is \(\kappa_f\ne0\), so it is a diffeomorphism onto its image.
All cutoffs below are supported in this interior, with a fixed
buffer around \(K\).

Pham's local fold calculation now applies: a degenerate point of the
cross-partition canonical relation satisfies equality of both \(f\)
and \(\rho\) at the two source points.  The injectivity just proved
forces those points to coincide, and the nonvanishing derivative in
\eqref{eq:Pham-level-monotonicity} supplies the fold transversality.
Thus the canonical relation corresponding to the cross
partition
\[
 (x,y';x',y)
\]
is then a two-sided Whitney fold; see
\cite[Lemma~4.2]{PhamFIO}.  In the notation of
\cite[Theorem~3.1]{PhamFIO}, the target dimension is $p=1$, the
partition parameter is $\alpha=2$, and the fold loses
$\beta=1/6$ derivatives.  The resulting total loss is
\begin{equation}
 \frac{\alpha-p}{2}+\beta=\frac{1}{2}+\frac{1}{6}=\frac{2}{3}.
 \label{eq_pham_two_thirds_loss}
\end{equation}

We include the Sobolev calculation because the norm inequality
\eqref{eq_pham_local_energy} is stronger than the support-dimension
statement of Pham's theorem.  Put
\[
 \sigma=\frac{a_1+a_2-2}{2},
 \qquad \mu=\mu_1\otimes\mu_2.
\]
The Fourier characterization of energy gives
\begin{equation}
 \|\mu\|_{H^\sigma(\mathbb R^2)}^2
 \ll
 \bigl(1+I_{a_1}(\mu_1)\bigr)
 \bigl(1+I_{a_2}(\mu_2)\bigr).
 \label{eq_pham_product_sobolev}
\end{equation}
For example, this follows from
\[
 (1+|\xi|^2+|\eta|^2)^\sigma
 \ll (1+|\xi|^2)^{(a_1-1)/2}
      (1+|\eta|^2)^{(a_2-1)/2}.
\]
For the output weight $(1+|\theta|)^{u-1}$, the effective operator
order in Pham's proof is $m_{\rm eff}=u-3/2$.
{Indeed, the amplitude has order \(u-1\), there is one phase
variable, and the left and right spatial dimensions are both \(2\).
Thus the H\"ormander order is
\[
 (u-1)+\frac{1}{2}-\frac{2+2}{4}=u-\frac{3}{2}.
\]
Since the partition is balanced, its effective order is the same.}
Consequently, the Sobolev pairing condition is
\begin{equation}
 2\sigma-m_{\rm eff}-\frac{1}{6}
 =a_1+a_2-u-\frac{2}{3}>0.
 \label{eq_pham_pairing_condition}
\end{equation}

We now define the operator used in the cross-partition step.  Choose
a smooth cutoff \(\chi(x,y,x',y')\), supported where both
\((x,y)\) and \((x',y')\) lie in \(K^+\), and equal to one when both
points lie in \(K\).  With left variables \((x,y')\) and right
variables \((x',y)\), define
\begin{align}
 (Tg)(x,y')
 :=\iiint &e^{-2\pi i\theta(f(x,y)-f(x',y'))}
 \langle\theta\rangle^{u-1}
 \chi(x,y,x',y')\notag\\
 &\hspace{25mm}\cdot g(x',y)\,d\theta\,dx'\,dy.
 \label{eq_pham_cross_operator}
\end{align}
where \(\langle\theta\rangle=(1+|\theta|^2)^{1/2}\).
The integral is first understood with a smooth frequency truncation.
Set
\[
 \mu_L=\mu_1(x)\otimes\mu_2(y'),
 \qquad
 \mu_R=\mu_1(x')\otimes\mu_2(y).
\]
Fubini's theorem for smooth input densities gives the exact pairing
identity
\begin{equation}
 \big\langle T\mu_R,\mu_L\big\rangle
 =\int_{\mathbb R}|\widehat\nu(\theta)|^2
 \langle\theta\rangle^{u-1}\,d\theta,
 \qquad \nu=f_*(\mu_1\times\mu_2).
 \label{eq_pham_cross_partition_pairing}
\end{equation}

The proof of \cite[Theorem~3.1]{PhamFIO}, in particular the
calculation in \cite[(3.8)--(3.9)]{PhamFIO}, decomposes the operator
\eqref{eq_pham_cross_operator} into finitely many microlocal pieces.
By \cite[Lemma~4.2]{PhamFIO}, every relevant canonical relation is a
two-sided Whitney fold.  The fold estimate
\cite[Theorem~2.2]{PhamFIO} therefore gives, for each piece and hence
for \(T\),
\begin{equation}
 T:H^\sigma_{\rm comp}(\mathbb R^2)
 \longrightarrow
 H^{\sigma-m_{\rm eff}-1/6}_{\rm loc}(\mathbb R^2),
 \label{eq_pham_cross_mapping}
\end{equation}
with a norm depending only on the fixed cutoffs and the buffered
chart.  Condition \eqref{eq_pham_pairing_condition} implies the
continuous embedding
\[
 H^{\sigma-m_{\rm eff}-1/6}\hookrightarrow H^{-\sigma}.
\]
Sobolev duality, \eqref{eq_pham_cross_partition_pairing}, and
\eqref{eq_pham_product_sobolev} now give
\begin{align*}
 \int_{\mathbb R}|\widehat\nu(\theta)|^2
 \langle\theta\rangle^{u-1}\,d\theta
 &\lesssim
 \|\mu_R\|_{H^\sigma}\|\mu_L\|_{H^\sigma}\\
 &\lesssim
 \bigl(1+I_{a_1}(\mu_1)\bigr)
 \bigl(1+I_{a_2}(\mu_2)\bigr).
\end{align*}
Here \(\mu_L\) and \(\mu_R\) have the same product Sobolev norm as
\(\mu\), up to the harmless permutation of variables.
For compactly supported measures and $0<u<1$, the last integral is
equivalent, up to a bounded low-frequency term, to $I_u(\nu)$.
This proves \eqref{eq_pham_local_energy}.

The preceding argument may first be applied to smooth measures
supported in the interiors of $I$ and $J$.  Approximation by such
measures and lower semicontinuity of Riesz energy give the stated
inequality for arbitrary finite-energy probability measures.  If an
input energy is infinite, \eqref{eq_pham_local_energy} is immediate.

We emphasize the source of the assertion: the norm estimate
\eqref{eq_pham_local_energy} is not stated verbatim in
\cite{PhamFIO}.  It is obtained by keeping the operator norms in
Pham's proof of Theorem~3.1.  The dimension conclusion of that theorem
alone would not imply \eqref{eq_pham_local_energy}.
\end{proof}

We now pass from independent product measures to a possibly dependent
coupling.  This is the main reduction used in both proofs.

\begin{proposition}[Correlation penalty]
\label{prop_pham_correlation}
Let $f\in\mathbb R[x,y]$ be non-special.  Suppose that $(U,V)$ has
bounded range and is jointly $(a,b,K)$-Frostman, where
\[
 0<a,b\leq1,\qquad a+b>\frac{4}{3}.
\]
Put
\begin{equation}
 q(a,b)=\min\left\{1,a+b-\frac{2}{3}\right\}.
 \label{eq_pham_q_ab}
\end{equation}
For every $0<t<q(a,b)$, there are a fixed integer $r\geq0$ and a
constant $B<\infty$ such that
\begin{equation}
 H_n(f(U,V))+2\operatorname{I}_{n+r}(U;V)\geq t n-B
 \label{eq_pham_correlation}
\end{equation}
for every $n\geq1$.  The constants are uniform over all couplings
with the fixed displayed data.
\end{proposition}

\begin{proof}
Choose
\[
 t<u<q(a,b),\qquad a'<a,\quad b'<b,
 \qquad a'+b'>u+\frac{2}{3}.
\]
We first assume that the law $\mathsf P$ of $(U,V)$ is supported on one
product rectangle compactly contained in
$\mathbb R^2\setminus Z(P_f)$.  Let
$\mathsf Q=\mathsf P_U\times\mathsf P_V$.  Quantize the
two source coordinates at level $n+r$, where the fixed integer $r$
is chosen so that replacing each source cell by a representative
changes $f$ by $O(2^{-n})$.  Let $p$ and $q$ be the resulting
level-$n$ output laws under $\mathsf P$ and $\mathsf Q$, respectively.

The marginal laws are $a$- and $b$-Frostman.  Hence
\eqref{eq_pham_input_energy},
Proposition~\ref{prop_pham_local_energy}, and
\eqref{eq_pham_energy_entropy} give
\begin{equation}
 H_{\mathrm{coll}}(q)\geq u n-O(1).
 \label{eq_pham_product_collision}
\end{equation}
Replacing representatives by the original variables changes Shannon
entropy by $O(1)$.  It changes collision entropy by at most $O(1)$ as
well: bounded displacement ensures both that every old bin reaches
only $O(1)$ new bins and that every new bin receives mass from only
$O(1)$ old bins.  Cauchy--Schwarz therefore changes the squared
$\ell^2$ norm by at most a fixed factor.

Data processing for relative entropy gives
\begin{equation}
 D(p\Vert q)
 \leq D\big(\mathsf P_{n+r}\Vert
       (\mathsf P_U)_{n+r}\times(\mathsf P_V)_{n+r}\big)
 =\operatorname{I}_{n+r}(U;V).
 \label{eq_pham_data_processing}
\end{equation}
Lemma~\ref{lem_pham_relative_entropy},
\eqref{eq_pham_product_collision}, and the bounded displacement just
described imply
\begin{equation}
 H_n(f(U,V))+2\operatorname{I}_{n+r}(U;V)
 \geq u n-O(1)
 \label{eq_pham_local_correlation}
\end{equation}
on this rectangle.

It remains to localize uniformly.  A fixed real algebraic curve in a
bounded square has an $\eta$-neighborhood that can be covered by
$O(\eta^{-1})$ squares of side $O(\eta)$.  It follows from the joint
Frostman condition that
\begin{equation}
 \mathbb P\big(\operatorname{dist}((U,V),Z(P_f))<\eta\big)
 \ll \eta^{a+b-1}.
 \label{eq_pham_bad_tube}
\end{equation}
Fix $0<\delta<1-t/u$.  Choose $\eta>0$ so that the left-hand side of
\eqref{eq_pham_bad_tube} is at most $\delta$.  Partition the support
square into finitely many sufficiently small product cells.  Merge
the cells meeting the smaller bad neighborhood into one bad label;
each remaining cell has a fixed buffered enlargement disjoint from
$Z(P_f)$.  Let $L$ be the resulting finite label, let
$w_j=\mathbb P(L=j)$,
and call a label good when its cell is disjoint from the bad
neighborhood.  The total good mass $W$ satisfies
\begin{equation}
 W\geq1-\delta.
 \label{eq_pham_good_mass}
\end{equation}

Conditioning on a good cell of mass $w_j$ changes the two marginal
Frostman constants by at most a factor $O(1/w_j)$.  Applying the local
argument and \eqref{eq_pham_input_energy} therefore gives
\begin{align}
 &H(D_nf(U,V)\mid L=j)
 +2\operatorname{I}(D_{n+r}U;D_{n+r}V\mid L=j)\notag\\
 &\hspace{35mm}\geq u n-B_0-B_1\log\frac{1}{w_j},
 \label{eq_pham_cell_correlation}
\end{align}
where $B_0,B_1$ do not depend on $j$ or on the coupling.  The same
scale shift $r$ works for all good labels because there are only
finitely many cells.

Entropy concavity and conditional mutual information give
\begin{align*}
 \sum_{j\ \mathrm{good}}w_jH(D_nf(U,V)\mid L=j)
 &\leq H_n(f(U,V)),\\
 \sum_{j\ \mathrm{good}}w_j
 \operatorname{I}(D_{n+r}U;D_{n+r}V\mid L=j)
 &\leq\operatorname{I}(D_{n+r}U;D_{n+r}V\mid L)\\
 &\leq\operatorname{I}_{n+r}(U;V)+H(L).
\end{align*}
Moreover,
\[
 \sum_jw_j\log\frac{1}{w_j}=H(L)\leq\log|\operatorname{range}(L)|=O(1).
\]
Averaging \eqref{eq_pham_cell_correlation} over the good labels yields
\[
 H_n(f(U,V))+2\operatorname{I}_{n+r}(U;V)
 \geq Wu n-O(1)>t n-O(1),
\]
by the choice of $\delta$.  This proves
\eqref{eq_pham_correlation}.
\end{proof}

\subsection{Proof of the explicit expansion theorem}

\begin{proof}[Proof of \cref{thm:Pham-explicit}\textup{(i)}]
Set
\[
 q_0=\min\left\{1,2s-\frac{2}{3}\right\}.
\]
Fix \(t\) with \(0<t<q_0\), and let \(r\) be the scale shift in
\cref{prop_pham_correlation} for \(f\) and \(t\).  Write
\[
 M_n=\max\left\{\frac{H_n(X)+H_n(Y)}{2},H_n(f(X,Y))\right\},
 \qquad m=n+r.
\]
By \eqref{eq:entropy-refinement} and
\cref{lem:Frostman-entropy},
\begin{align*}
 \operatorname{I}_m(X;Y)
 &=H_m(X)+H_m(Y)-H_m(X,Y)\\
 &\le 2M_n-2sn+O_r(1).
\end{align*}
The correlation penalty therefore gives
\[
 t n-O(1)
 \le H_n(f(X,Y))+2\operatorname{I}_m(X;Y)
 \le 5M_n-4sn+O_r(1).
\]
Consequently
\[
 M_n\ge\left(s+\frac{t-s}{5}\right)n-O(1).
\]
Since
\[
 q_0-s=\min\left\{s-\frac{2}{3},1-s\right\},
\]
we may choose \(s+5\varepsilon<t<q_0\) for every \(\varepsilon\) in the range
stated in \cref{thm:Pham-explicit}\textup{(i)}.  This proves that part.
\end{proof}

\begin{proof}[Proof of \cref{thm:Pham-explicit}\textup{(ii)}]
If \(f\) is not analytic-special, part \textup{(i)} applies.  Its
admissible range contains the full range in part \textup{(ii)}, since
\[
 \frac{1}{7}\min\left\{2s-\frac{5}{3},1-s\right\}
 \le
 \frac{1}{5}\min\left\{s-\frac{2}{3},1-s\right\}.
\]
We may therefore assume that \(f\) is analytic-special.  Since
\(f\notin\mathfrak E_3\), \cref{cor:three-coordinate} supplies a
non-special sum transition.  Thus, with \(S=X+Y\), there are
\[
 (U,V)=(X,Y)\quad\hbox{or}\quad(U,V)=(Y,X)
\]
and a non-special \(G\) such that \(G(U,S)=f(X,Y)\).
By \cref{lem:joint-shear}, \((U,S)\) is jointly
\((s,2s-1,C')\)-Frostman.  The limiting exponent in
\cref{prop_pham_correlation} is
\[
 q_1=\min\left\{1,3s-\frac{5}{3}\right\},
 \qquad
 q_1-s=\min\left\{2s-\frac{5}{3},1-s\right\}.
\]
The assumption \(s>5/6\) ensures both that \(q_1>s\) and that
\(s+(2s-1)>4/3\).

Let
\[
 M_n=\max\left\{\frac{H_n(X)+H_n(Y)}{2},H_n(S),
                         H_n(f(X,Y))\right\}.
\]
Fix \(0<t<q_1\), let \(r\) be the shift in the correlation penalty for
\(G,(U,S)\), and put \(m=n+r\).  If \(V\) denotes the unused original
coordinate, then
\[
 H_m(U)+H_m(V)\le2M_n+O_r(1),
 \qquad H_m(V)\ge sm-O(1).
\]
Moreover, the shear is invertible, so
\[
 H_m(U,S)=H_m(X,Y)+O(1)\ge2sm-O(1).
\]
Together with \(H_m(S)\le M_n+O_r(1)\), these inequalities give
\[
 \operatorname{I}_m(U;S)\le3M_n-3sn+O_r(1).
\]
Applying \cref{prop_pham_correlation} to \(G(U,S)=f(X,Y)\) yields
\[
 t n-O(1)\le H_n(f(X,Y))+2\operatorname{I}_m(U;S)
          \le7M_n-6sn+O_r(1).
\]
Thus
\[
 M_n\ge\left(s+\frac{t-s}{7}\right)n-O(1).
\]
Letting \(t\) approach \(q_1\) proves every gain strictly smaller
than \((q_1-s)/7\); in particular the safe value
\eqref{eq:Pham-joint-eps} follows.

Suppose now that both sum transitions are non-special.  Apply the
correlation penalty to both \((X,S)\) and \((Y,S)\), taking a common
larger fixed scale shift \(r\), which is allowed by
\eqref{eq:Pham-information-refinement}.  The invertible-shear
estimate and the preceding bounds give
\begin{align*}
 \operatorname{I}_m(X;S)+\operatorname{I}_m(Y;S)
 &=H_m(X)+H_m(Y)+2H_m(S)\\
 &\quad-H_m(X,S)-H_m(Y,S)\\
 &\le4M_n-4sn+O_r(1).
\end{align*}
One of these two mutual informations is therefore at most
\(2M_n-2sn+O_r(1)\).  Using the corresponding transition in the
correlation penalty gives
\[
 t n-O(1)\le5M_n-4s n+O_r(1).
\]
Hence every gain smaller than \((q_1-s)/5\) is valid, and the safe
factor \(1/10\) asserted in part \textup{(ii)} follows.
\end{proof}

\begin{proof}[Proof of \cref{thm:Pham-explicit}\textup{(iii)}]
Put \(S=X+Y\), and write
\[
 M_n=\max\{H_n(S),H_n(f(X,Y))\},\qquad
 \bar s=\frac{s_1+s_2}{2},\qquad
 q=\min\left\{1,s_1+s_2-\frac{2}{3}\right\}.
\]
The hypotheses give \(2/3<\bar s<1\), and hence
\[
 q-\bar s=\min\left\{\bar s-\frac{2}{3},1-\bar s\right\}>0.
\]
The conditional-shear estimates in \cref{lem:conditional-shear}
also give the baseline
\begin{equation}
 H_n(S)\ge\max\{s_1,s_2\}n-O(1).
 \label{eq:Pham-sum-asymmetric-baseline}
\end{equation}

First suppose that \(f\notin\mathfrak E_S\).  By
\cref{thm:two-transition}, at least one sum transition is non-special.
Write it as \(G(U,S)=f(X,Y)\), with
\[
 (U,V;s_U,s_V)=(X,Y;s_1,s_2)
 \quad\hbox{or}\quad(Y,X;s_2,s_1).
\]
By \cref{lem:conditional-shear}, \((U,S)\) is jointly
\((s_U,s_V,C')\)-Frostman, and
\begin{equation}
 H(D_mS\mid D_mU)\ge s_Vm-O(1).
 \label{eq:Pham-transition-conditional-entropy}
\end{equation}
Fix \(0<t<q\), and put \(m=n+r\), where \(r\) is the fixed
shift in \cref{prop_pham_correlation} for \(G,(U,S)\).
Then
\[
 \operatorname{I}_m(U;S)
 =H_m(S)-H(D_mS\mid D_mU)
 \le M_n-s_Vn+O_r(1).
\]
The correlation penalty implies
\begin{equation}
 3M_n\ge(t+2s_V)n-O(1).
 \label{eq:Pham-oriented-three-bound}
\end{equation}
Combining this with \(M_n\ge s_Un-O(1)\), which follows from
\eqref{eq:Pham-sum-asymmetric-baseline}, proves
\eqref{eq:Pham-cond-oriented} and its symmetric counterpart.
Adding twice the latter baseline to
\eqref{eq:Pham-oriented-three-bound} also yields
\begin{equation}
 5M_n\ge\bigl(t+2(s_1+s_2)\bigr)n-O(1).
 \label{eq:Pham-asymmetric-five-bound}
\end{equation}
Therefore
\[
 M_n\ge\left(\bar s+\frac{t-\bar s}{5}\right)n-O(1).
\]
For every \(0<\varepsilon<(q-\bar s)/5\), choose
\(\bar s+5\varepsilon<t<q\).  This proves
\eqref{eq:Pham-cond-main} in the transition case.

The transition-specific bound is retained as well.  Since
\(s_V\ge s_-\), \eqref{eq:Pham-oriented-three-bound} gives
\[
 M_n\ge\left(s_-+\frac{t-s_-}{3}\right)n-O(1).
\]
Choosing \(s_-+3\varepsilon<t<q\) proves
\eqref{eq:Pham-cond-transition-main} for every
\(0<\varepsilon<(q-s_-)/3\).

It remains to treat \(f\in\mathfrak E_S\setminus\mathfrak E_3\).
Both sum transitions are analytic-special, so
\cref{cor:three-coordinate} shows that \(f\) itself is non-special.
Moreover \(f_y-f_x\not\equiv0\), since otherwise
\(f=Q(x+y)\in\mathfrak E_3\).  Fix \(0<t<q\), and let
\(m=n+r\) be the scale in \cref{prop_pham_correlation} for
\(f,(X,Y)\).  By \cref{thm:coordinate-entropy} and
\eqref{eq:entropy-refinement},
\[
 H_m(X,Y)=H_m(S,f(X,Y))+O(1)\le2M_n+O_r(1).
\]
The conditional entropy bounds give
\begin{align*}
 \operatorname{I}_m(X;Y)
 &=H_m(X,Y)-H(D_mX\mid D_mY)-H(D_mY\mid D_mX)\\
 &\le2M_n-(s_1+s_2)n+O_r(1).
\end{align*}
Applying the correlation penalty directly to \(f(X,Y)\), we obtain
\[
 tn-O(1)\le H_n(f(X,Y))+2\operatorname{I}_m(X;Y)
 \le5M_n-2(s_1+s_2)n+O_r(1).
\]
This is again \eqref{eq:Pham-asymmetric-five-bound}, and choosing
\(\bar s+5\varepsilon<t<q\) proves the full range in
\eqref{eq:Pham-cond-main}.  Taking half of each open gain range
gives the stated safe values.  When \(s_1=s_2=s\), one has
\(\bar s=s_-=s\); the uniform and transition-specific factors are
therefore \(1/10\) and \(1/6\), respectively.
\end{proof}

\begin{proof}[Proof of \cref{thm:Pham-explicit}\textup{(iv)}]
First suppose that \(f\) is not analytic-special.  The assertion is
vacuous if
\(\min\{1,s_1+s_2-2/3\}\le0\), so fix
\[
 0<u<u'<\min\left\{1,s_1+s_2-\frac{2}{3}\right\}.
\]
Choose \(a_i<s_i\) so close to \(s_i\) that
\(a_1+a_2>u'+2/3\).

The independent pair \((X,Y)\) is conditionally
\((s_1,C_1;s_2,C_2)\)-Frostman.  Apply
\cref{lem:algebraic-localization} to the nonzero polynomial \(P_f\),
with bad mass at most \(\beta\), where \(\beta>0\) is chosen so that
\((1-\beta)u'>u\).  Each good label is a product rectangle.
Independence is preserved after conditioning on such a rectangle,
so the conditional law is the product of its two conditional
marginals.  The Frostman constants of those marginals grow by at
most the reciprocal of their conditioning masses.  Thus
\eqref{eq_pham_input_energy},
\cref{prop_pham_local_energy}, and
\eqref{eq_pham_energy_entropy} give, on every good label \(j\),
\[
 H_n(f(X,Y)\mid L=j)
 \ge u'n-B_0-B_1\log\frac{1}{w_j},
 \qquad w_j=\mathbb P(L=j),
\]
with constants independent of \(j,n\).  Averaging by
\eqref{eq:mixture-concavity} and using
\(\sum_jw_j\log(1/w_j)\le H(L)=O(1)\), we get
\[
 H_n(f(X,Y))\ge(1-\beta)u'n-O(1)\ge u n-O(1).
\]
This proves the first assertion of part \textup{(iv)}.

Finally suppose only that \(f\notin\mathfrak E_3\) and
\(s_1+s_2>4/3\).  Independence and the marginal Frostman bounds
imply the conditional hypotheses of part \textup{(iii)}, so its
conclusions apply, including the bound above \(\bar s\) and the
additional transition-specific estimates when their hypotheses hold.
\end{proof}

\section{Examples and genuine obstructions}
\label{sec:examples-obstructions}

We conclude with examples that illustrate the scope of
\cref{thm:RZ-general,thm:Pham-explicit}.  The first group places
familiar diagonal and quadratic polynomials inside the general
theorems and exhibits the two transition alternatives in
\cref{thm:Pham-explicit}.  The second example explains why the
original coordinates must be used together with the two sum-coordinate
systems.  The final proposition shows that the excluded family
\(\mathfrak E_3\) contains genuine obstructions, even when the
Frostman constant is fixed.

\begin{proposition}[Diagonal and quadratic examples]
\label{prop:diagonal-quadratic-examples}
The following polynomials lie outside \(\mathfrak E_3\).
\begin{enumerate}[label=\textup{(\roman*)}]
\item
\begin{equation}
 f(x,y)=\rho_1(x)+\rho_2(y),
 \label{eq:diagonal-example}
\end{equation}
where \(\rho_1,\rho_2\in\mathbb R[t]\) are nonconstant and at least one of
them has degree at least two.  In particular, this includes
\(x^d+y^d\) for every \(d\ge2\).

\item
\begin{equation}
 f(x,y)=Ax^2+2Bxy+Cy^2+Dx+Ey+F,
 \label{eq:quadratic-example}
\end{equation}
provided that at least one of
\begin{equation}
 AC-B^2,\qquad AE-BD,\qquad BE-CD
 \label{eq:quadratic-example-minors}
\end{equation}
is nonzero.  In particular, every quadratic polynomial whose
quadratic part is nondegenerate belongs to this class.
\end{enumerate}
Consequently, the qualitative conclusions of
\cref{thm:RZ-general} apply to these polynomials throughout
\(0<s<1\).  The conclusions of
\cref{thm:Pham-explicit}\textup{(ii)--(iv)} that assume
\(f\notin\mathfrak E_3\) also apply in their stated Frostman ranges.
\end{proposition}

\begin{proof}
For \textup{(i)}, suppose that
\(u f_x+v f_y\equiv0\).  Then
\[
 u\rho_1'(x)+v\rho_2'(y)\equiv0.
\]
Differentiating in a variable for which the corresponding polynomial
is nonlinear shows that its coefficient, \(u\) or \(v\), vanishes.
The original identity and nonconstancy of the other polynomial then
show that the other coefficient also vanishes.  Thus
\cref{cor:E3-test} gives \(f\notin\mathfrak E_3\).

For \textup{(ii)}, the identity \(u f_x+v f_y\equiv0\) is equivalent
to
\[
 \begin{pmatrix}
 A&B\\ B&C\\ D&E
 \end{pmatrix}
 \begin{pmatrix}u\\v\end{pmatrix}=0.
\]
A nonzero solution exists precisely when this \(3\times2\) matrix
has rank at most one, which is equivalent to the vanishing of all
three minors in \eqref{eq:quadratic-example-minors}.  The conclusion
again follows from \cref{cor:E3-test}.
\end{proof}

The two alternatives in the jointly Frostman transition argument
already occur among diagonal quadratics.  For
\[
 f_0(x,y)=x^2+y,
\]
the sum transitions are
\[
 G_{0,X}(x,S)=x^2-x+S,
 \qquad
 G_{0,Y}(y,S)=(S-y)^2+y.
\]
The first is additively special, whereas a direct calculation from
\eqref{eq_pham_cleared_curvature} gives
\[
 \mathcal B_{G_{0,Y}}=16(S-y)-4\not\equiv0.
\]
Thus exactly one transition is non-special.  This is the default
case of \cref{thm:Pham-explicit}\textup{(ii)}, with the safe factor
\(1/14\); by \cref{thm:two-transition},
\(f_0\notin\mathfrak E_S\), so it also falls under the
stronger range in \cref{thm:Pham-explicit}\textup{(iii)}.  In
contrast, for
\[
 f_1(x,y)=x^2+y^2,
\]
the two transitions are
\[
 G_{1,X}(x,S)=2x^2-2xS+S^2,
 \qquad
 G_{1,Y}(y,S)=2y^2-2yS+S^2,
\]
and
\[
 \mathcal B_{G_{1,X}}=16(S^2-2x^2),
 \qquad
 \mathcal B_{G_{1,Y}}=16(S^2-2y^2).
\]
Both transitions are non-special, so the improved safe factor
\(1/10\) in \cref{thm:Pham-explicit}\textup{(ii)} applies.

\begin{example}[Why the original coordinates are needed]
\label{ex:original-coordinates-needed}
Let
\[
 f(x,y)=(x+y)^2+x-y.
\]
Writing \(S=x+y\), we have
\[
 f(x,y)=2x+S^2-S,
\]
so \(f\in\mathfrak E_S\).  Its two sum transitions are
\[
 G_X(x,S)=S^2+2x-S,
 \qquad
 G_Y(y,S)=S^2+S-2y,
\]
and hence both are additively special.  On the other hand,
\[
 f_x=2(x+y)+1,\qquad f_y=2(x+y)-1.
\]
If \(u f_x+v f_y\equiv0\), comparison of the coefficient of
\(x+y\) and the constant term gives \(u+v=u-v=0\).  Therefore
\(f\notin\mathfrak E_3\).  Moreover,
\begin{equation}
 \mathcal B_f=32(x+y)\not\equiv0,
 \label{eq:original-coordinate-curvature-example}
\end{equation}
so the original representation is non-special.

Thus \cref{thm:RZ-general} applies for every \(0<s<1\), even though
both sum transitions are special.  In the jointly Frostman setting,
\cref{thm:Pham-explicit}\textup{(i)} applies for \(2/3<s<1\), with
every
\[
 0<\varepsilon<\frac{1}{5}\min\left\{s-\frac{2}{3},1-s\right\}.
\]
This example shows concretely why the three-coordinate exceptional
family is \(\mathfrak E_3\), rather than the larger family
\(\mathfrak E_S\).
\end{example}

\begin{proposition}[Algebraic affine-univariate obstructions]
\label{prop:algebraic-affine-obstructions}
Let
\[
 f(x,y)=\rho(ax+by+c_0),
 \qquad (a,b)\ne(0,0),
\]
where \(\rho\in\mathbb R[t]\), and suppose that the projective direction
\([a:b]\) is algebraic; equivalently, either \(a=0\), or
\(b/a\) is a real algebraic number.  Fix \(0<s<1\) and \(c>0\).
There are a constant \(C_0=C_0(f,s,c)\), arbitrarily large integers
\(n\), and independent \((s,C_0)\)-Frostman random variables
\(X_n,Y_n\), supported in \([-c,c]\), such that
\begin{equation}
 \max\left\{
 \frac{H_n(X_n)+H_n(Y_n)}{2},\,
 H_n(X_n+Y_n),\,
 H_n(f(X_n,Y_n))
 \right\}
 \le sn+O_{f,s,c}(1).
 \label{eq:algebraic-affine-obstruction}
\end{equation}
Consequently, the hypotheses excluding \(\mathfrak E_3\) in
\cref{thm:RZ-general,thm:Pham-explicit} cannot simply be removed.
\end{proposition}

\begin{proof}
After interchanging \(X_n\) and \(Y_n\), and simultaneously
interchanging \(a\) and \(b\), if necessary, we may assume that
\(a\ne0\).  Put \(\lambda=b/a\).  This includes the case \(b=0\),
for which \(\lambda=0\), \(K=\mathbb Q\), and \(d=1\) below.  Let
\(K=\mathbb Q(\lambda)\), let \(d=[K:\mathbb Q]\), and choose an
integral basis \(\omega_1,\ldots,\omega_d\) of the ring of integers
of \(K\).  For \(N\ge2\), set
\begin{equation}
 A_N=\left\{\frac{1}{K_0N}\sum_{j=1}^d m_j\omega_j:
 0\le m_j<N\right\},
 \qquad M=N^d,
 \label{eq:algebraic-module-box}
\end{equation}
where the fixed normalization \(K_0\) is chosen so that, for all
large \(N\), the small neighborhoods of \(A_N\) used below are
contained in \([-c,c]\).

The elements of \(A_N\) are distinct.  If
\(\alpha=\sum_jk_j\omega_j\ne0\), with \(|k_j|<N\), then its
algebraic norm is a nonzero integer.  Every conjugate of \(\alpha\)
other than the distinguished real conjugate has absolute value
\(O_K(N)\).  Hence
\begin{equation}
 |\alpha|\gtrsim_K N^{-(d-1)},
 \qquad
 \operatorname{sep}(A_N)\gtrsim_{K,c}N^{-d}=M^{-1}.
 \label{eq:algebraic-module-separation}
\end{equation}

Put \(\delta=M^{-1/s}\), and let \(\mu_N\) be the convolution of the
uniform probability measure on \(A_N\) with normalized Lebesgue
measure on an interval of length \(c_1\delta\), where \(c_1>0\) is
fixed and sufficiently small.  Since \(s<1\), these intervals are
disjoint for all sufficiently large \(N\).  If \(I\) is an interval
with \(\delta\le |I|\le1\), then
\eqref{eq:algebraic-module-separation} gives
\[
 \mu_N(I)\lesssim |I|+M^{-1}\lesssim |I|^s.
\]
If \(|I|<\delta\), the density on each component is
\(O((M\delta)^{-1})\), and therefore
\[
 \mu_N(I)\lesssim\frac{|I|}{M\delta}\le |I|^s,
\]
because
\(
 |I|^{1-s}\le\delta^{1-s}=M\delta
\).
Thus \(\mu_N\) is
\((s,C_0)\)-Frostman with \(C_0\) independent of \(N\).

Addition is a fixed linear operation on the coefficient box in
\eqref{eq:algebraic-module-box}.  Multiplication by \(\lambda\) is a
fixed rational linear map on the coefficient lattice with respect to
the basis \(\omega_1,\ldots,\omega_d\).  Clearing its fixed
denominators and counting coefficient vectors gives
\begin{equation}
 |A_N+A_N|+|A_N+\lambda A_N|\lesssim_\lambda N^d
 =O_\lambda(M).
 \label{eq:algebraic-module-small-images}
\end{equation}
Let \(X_n,Y_n\) be independent with law \(\mu_N\), and choose the
integer \(n\) so that \(2^{-n}\asymp\delta\).  Each smoothing
component meets only \(O(1)\) level-\(n\) intervals.  It follows from
\eqref{eq:algebraic-module-small-images} that
\begin{align*}
 H_n(X_n),H_n(Y_n),H_n(X_n+Y_n),
 H_n(X_n+\lambda Y_n)
 &\le\log M+O_{K,s,c}(1)\\
 &=sn+O_{K,s,c}(1).
\end{align*}
Finally,
\[
 f(X_n,Y_n)=\rho\bigl(a(X_n+\lambda Y_n)+c_0\bigr),
\]
and \(\rho\) is Lipschitz on the fixed relevant interval.  Hence
\eqref{eq:Lipschitz-entropy} gives
\eqref{eq:algebraic-affine-obstruction}.

Independence makes the pair conditionally
\((s,C_0;s,C_0)\)-Frostman.  It is also jointly
\((s,s,C_0^2)\)-Frostman.  Thus the same construction obstructs the
corresponding conclusions in
\cref{thm:RZ-general,thm:Pham-explicit}.
\end{proof}

The preceding proposition gives genuine obstructions for every
algebraic projective direction.  We do not claim here that every
element of \(\mathfrak E_3\) is a counterexample; the point is that
the exceptional family in the main theorems is not merely a formal
artifact of the proofs.

\appendix
\section{An alternative proof using curved three-webs}
\label{sec:curved-appendix}

For completeness, we give a second proof of the non-special branch
of \cref{thm:RZ-general}
which uses the curved three-web theorem in place of
Proposition~\ref{prop:RZ-global}.  This proof is available under the
conditional Frostman hypothesis for every \(0<s<1\), and under the
joint Frostman hypothesis when \(s>1/2\).

Assume that \(f\) is not a polynomial special form and that one of
the conclusions of \cref{thm:RZ-general} fails at level \(n\).
Use the normalization, entropy bounds, and information profiles from
\cref{sec:RZ-proof}, with \(\delta=2^{-n}\).  In particular,
\eqref{eq:RZ-profile-size}--\eqref{eq:RZ-profile-image-upper}
are available before any chart labels are adjoined.

Set
\begin{equation}
 \mathcal K_f
 =f_y^2\bigl(f_xf_{xxy}-f_{xx}f_{xy}\bigr)
  -f_x^2\bigl(f_yf_{xyy}-f_{xy}f_{yy}\bigr).
 \label{eq:RZ-curvature-numerator}
\end{equation}
On the set where \(f_xf_y\ne0\), the Blaschke curvature of the
three-web \((x,y,f)\) is, up to a fixed nonzero normalization,
\begin{equation}
 \partial_{xy}\ln\left|\frac{f_x}{f_y}\right|
 =\frac{\mathcal K_f}{f_x^2f_y^2}.
 \label{eq:RZ-curvature-identity}
\end{equation}
If \(f\) is not a polynomial special form, then none of
\(f_x,f_y,f_{xy}\) vanishes identically, and the polynomial form of
the Elekes--R\'onyai curvature criterion
\cite[Lemma~3.3]{RZ} shows that
\(\mathcal K_f\not\equiv0\).  Thus
\[
 R_f=f_xf_y\mathcal K_f
\]
is a nonzero polynomial.

Apply Lemma~\ref{lem:algebraic-localization} to \(R_f\), with bad
mass at most \(\beta\), and subdivide the resulting good rectangles
a fixed number of times.  On every buffered rectangle obtained in
this way, the gradients of \(x,y,f\) are pairwise transverse, the
curvature in \eqref{eq:RZ-curvature-identity} is bounded away from
zero, and each of the maps
\[
 (x,y),\qquad (x,f),\qquad (y,f)
\]
is bi-Lipschitz onto its image.  The localization lemma applies here
to conditionally Frostman laws for every \(s>0\), and to jointly
Frostman laws when \(s>1/2\).

Adjoin the finite rectangle label to the profiles from
\cref{sec:RZ-proof} and retain good children of mass at least
\(\delta^{2\zeta}\).  By \cref{rem:profile-refinement}, the
additional small-child loss is \(o(1)\), and the estimates
\eqref{eq:RZ-profile-size}, \eqref{eq:RZ-profile-source-nc}, and
\eqref{eq:RZ-profile-image-upper} remain valid for each child's
saturated source set \(E\).  Define
\[
 b_i=\frac{\log\mathcal N_\delta(\pi_i(E))}{\log(1/\delta)},
 \qquad i\in\{x,y,f\}.
\]
The three bi-Lipschitz coordinate pairs give
\begin{equation}
 b_x+b_y,\quad b_x+b_f,\quad b_y+b_f
 \ge2s-O(\zeta).
 \label{eq:RZ-curved-pair-bounds}
\end{equation}

Use the uniform exponents in
Proposition~\ref{prop:RZ-curved} on all the finitely many charts, and
choose
\[
 2\zeta<\eta_{\rm curv},
 \qquad \gamma+2\zeta<\eta_{\rm curv}.
\]
The profile size and source nonconcentration estimates then meet the
hypotheses of that proposition.  With
\(\kappa_c=\kappa_{\rm curv}/2\), and after increasing \(n\) to
absorb fixed constants, it follows that
\begin{equation}
 \max\{b_x,b_y,b_f\}\ge s+\kappa_c.
 \label{eq:RZ-curved-max}
\end{equation}
Subtracting \(s\) from the three variables, the pair bounds in
\eqref{eq:RZ-curved-pair-bounds} and
\eqref{eq:RZ-curved-max} immediately imply
\begin{equation}
 b_x+b_y+b_f\ge3s+\kappa_c-O(\zeta).
 \label{eq:RZ-curved-LP}
\end{equation}
Indeed, if one of the three variables is at least
\(s+\kappa_c\), the sum of the other two is at least
\(2s-O(\zeta)\).

Let \(a_x,a_y,a_f\) be the corresponding information-window lower
endpoints.  By \eqref{eq:RZ-profile-image-upper},
\(b_i\le a_i+\zeta+O(1/n)\).  If \(W_{\rm curv}\) denotes the
total mass of the retained good child profiles, then
\begin{equation}
 W_{\rm curv}\ge
 1-O\!\left(\beta+\frac{\varepsilon}{\gamma}+\frac{1}{K}\right)-o(1).
 \label{eq:RZ-curved-retained-mass}
\end{equation}
Averaging \eqref{eq:RZ-curved-LP} therefore gives
\begin{equation}
 H_n(X)+H_n(Y)+H_n(f(X,Y))
 \ge W_{\rm curv}
       \bigl(3s+\kappa_c-O(\gamma+\zeta)\bigr)n-O(1).
 \label{eq:RZ-curved-averaged-LP}
\end{equation}
Failure of the conditional conclusion gives the upper bound
\((3s+5\varepsilon)n+O(1)\) for the left-hand side, by
\eqref{eq:RZ-marginal-upper}; failure of the joint conclusion gives
the upper bound \((3s+3\varepsilon)n\).  First choose
\(\beta,\gamma,\zeta,K^{-1}\) sufficiently small in terms of the
uniform curved-web gain, then choose \(\varepsilon\) sufficiently small,
and finally take \(n\) large.  Either upper bound contradicts
\eqref{eq:RZ-curved-averaged-LP}.  This proves the claimed localized
alternative.

\end{document}